\documentclass{article}
\usepackage{graphicx} 
\usepackage[margin=1in]{geometry}
\usepackage[utf8]{inputenc}
\usepackage{amsmath}
\usepackage{amsthm}
\usepackage{amsfonts}
\usepackage{amssymb}
\usepackage{ragged2e}
\usepackage{authblk}
\usepackage{cite}
\DeclareMathOperator{\sech}{sech}
\providecommand{\keywords}[1]
{
  \small	
  \textbf{\textit{Keywords: }} #1
}

\providecommand{\subjclass}[1]
{
  \small	
  \textbf{\textit{Mathematics Subject Classification: }} #1
}

\newtheorem{theorem}{Theorem}

\newtheorem{lemma}[theorem]{Lemma}
\newtheorem{proposition}[theorem]{Proposition}
\title{Analytic expressions pertaining to certain arithmetical functions}
\author{\emph{Aung Phone Maw}}
\date{October 2024}

\begin{document}

\maketitle
\begin{abstract}
    We demonstrate the general outlines of a method for obtaining analytic expressions for certain types of general arithmetical sums. In particular, analytical expressions for a general arithmetical sum whose terms are summed over either the positive integer solutions $(a,b)$ of the Diophantine equation $kb^2+da^2 = N$ or the positive integer solutions $(a,b)$ of the Diophantine equation $kb^2-da^2 = N$ are derived.
\end{abstract}
.\\
.
\section{Introduction}

This study should be understood as part of the quest towards finding the exact analytic expressions of arithmetic functions. The most well-known example (although not representative of the expressions that we shall proceed to describe in our study) would be the Rademacher's formula for the partition function $p(n)$ \cite{Rademacher1938OnTP}, which is arrived by the extensive use of the circle method. Here, we shall derive analytic expressions for sums, whose terms are summed over the positive integer solutions of certain Diophantine equations. To derive these expression, we shall take a different trajectory other than the circle method. Particularly, we shall exploit a certain invertibility property of the following formal series of partial fractions

\[ \label{parfrac-series} \tag{1.1}
\sum\limits_{n=1}^{\infty} \frac{f(n)}{n+z}.
\]

And we shall make extensive usage of the following Mittag-Leffler expansion formula, which also appears in Ramanujan's notebook \cite[p.~237]{RamaII}.  

\[ \label{mittag-leffler-coth} \tag{1.2}
\frac{1}{2x^2}-\frac{\pi \cosh(\theta x)}{2x\sinh(\pi x)} = \sum\limits_{k=1}^{\infty}\frac{(-1)^{k-1}\cos(k \theta)}{k^2+x^2},\hspace{0.12cm} |\theta| \leq \pi,
\]
And we shall also make use of some consequential formulae of the above expansion, such as 

\[ \label{Tformula} \tag{1.3}
\sum\limits_{n=1}^{\infty}\frac{1}{z^2+(n^2+x)^2}= \frac{\pi}{4}T_{x,z}- \frac{1}{2(x^2+z^2)},
\]

\[ \label{Vformula} \tag{1.4}
\sum\limits_{n=1}^{\infty}\frac{(-1)^n}{z^2+(n^2+x)^2}= \frac{\pi}{2}V_{x,z}- \frac{1}{2(x^2+z^2)},
\]

where 
\[ \label{Tdef} \tag{1.5}
T_{M,t} = \frac{v_{M,t} \sinh(2\pi u_{M,t}) + u_{M,t}\sin(2\pi v_{M,t})}{t\sqrt{M^2+t^2}(\sinh^2(\pi u_{M,t}) + \sin^2(\pi v_{M,t}))},
\]

\[ \label{Vdef} \tag{1.6}
V_{M,t} = \frac{v_{M,t}\sinh(\pi u_{M,t}) \cos(\pi v_{M,t}) + u_{M,t}\cosh(\pi u_{M,t})\sin(\pi v_{M,t})}{t\sqrt{M^2+t^2}(\sinh^2(\pi u_{M,t})+ \sin^2(\pi v_{M,t}))},
\]

\[ \label{uvdef} \tag{1.7}
u_{M,t} = \sqrt{\frac{\sqrt{M^2+t^2}+M}{2}}, \text{ and } v_{M,t} = \sqrt{\frac{\sqrt{M^2+t^2}-M}{2}} \text{(Note that } v_{M,t} = u_{-M,t} \text{)}. \] 

From our definition of $u$ and $v$, note that we also have $\sqrt{M+it} = u_{M,t}+iv_{M,t}$ and $\sqrt{M-it}= u_{M,t}-iv_{M,t}$.

Before proceeding to state our main results, we shall also define, in the same spirit as before, $G_{M,t,k}$ to be the following expression

\[ \label{Gdef} \tag{1.8}
\frac{(M^2u_{M,t}-t^2u_{M,t}-2Mtv_{M,t})\sin(\frac{2\pi v_{M,t}}{\sqrt{k}})+(2Mtu_{M,t}+M^2v_{M,t}-t^2v_{M,t})\sinh(\frac{2\pi u_{M,t}}{\sqrt{k}})}{(M^2+t^2)^2\sqrt{M^2+t^2}(\sinh^2(\frac{\pi u_{M,t}}{\sqrt{k}}))+\sin^2(\frac{\pi v_{M,t}}{\sqrt{k}}))},
\]

and we shall denote $G_{M,t,1}$ equivalently as $G_{M,t}$. With the preliminary notations described above, some of our main results are stated as follows. The following statement gives an analytic expression for the arithmetical sum $\sum\limits_{\substack{da^2+kb^2=N \\ (a,b) \in \mathbb{N}^2}} \frac{1}{b^4}$.

\vspace{0.5cm}

 \begin{proposition} For natural numbers $N, d$, and $k$, and for real value $t>0$, we have

\begin{flalign*} 
    &\frac{1}{k^2}\sum\limits_{\substack{da^2+kb^2=N \\ (a,b) \in \mathbb{N}^2}} \frac{1}{b^4} 
    \\ &= \left(\frac{ \pi^2}{12} - \frac{\pi \coth(\pi t)}{4} + \frac{\pi \coth(\pi t)}{48 t^2} \right)\epsilon_N + \frac{\pi^2}{3k(e^{2\pi t}-1)}{\sum\limits_{a=1}^{\infty} }^*\frac{1}{N-da^2}+ \frac{\pi - \coth(\pi t)}{2\pi}{\sum\limits_{a=1}^{\infty} }^*\frac{1}{(N-da^2)^2} &\\&+ \frac{(-1)^{N}(1-t^2)\coth(\pi t)}{\pi}{\sum\limits_{a=1}^{\infty} }^* \frac{1}{(N-da^2)^2} + \frac{(-1)^{N-1} \pi^3 \coth(\pi t)}{12kt} {\sum\limits_{a=1}^{\infty}}^* \frac{1}{\sinh(\frac{\pi (N-da^2)}{2t})}  &\\&+ \frac{(-1)^{N}\coth(\pi t)(t-\frac{1}{t})}{4} {\sum\limits_{a=1}^{\infty}}^* \frac{1}{(N-da^2)\sinh(\frac{\pi (N-da^2)}{2t})} + \frac{(-1)^{N}\pi \coth(\pi t)}{8} {\sum\limits_{a=1}^{\infty}}^* \frac{\cosh(\frac{\pi(N-da^2)}{2t})}{\sinh^2(\frac{\pi (N-da^2)}{2t})} &\\&- \frac{\pi^2\coth(\pi t)}{3k} \sum\limits_{r=1}^{\infty} \sum\limits_{a=1}^{\infty} \frac{(-1)^{r-1}e^{-2\pi t r}(N-da^2)}{4t^2r^2+(N-da^2)^2}  - t\coth(\pi t) \sum\limits_{r=1}^{\infty} \sum\limits_{a=1}^{\infty} \frac{(-1)^{r-1}re^{-2\pi t r}}{4t^2r^2+(N-da^2)^2} &\\& - \frac{\coth(\pi t)}{\pi}\sum\limits_{r=1}^{\infty}\sum\limits_{a=1}^{\infty} \frac{(-1)^{r-1}e^{-2\pi t r}(4t^2r^2-(N-da^2)^2)}{(4t^2r^2+(N-da^2)^2)^2} + \frac{\sinh(\pi t) G_{-N,t,k}}{8\sqrt{k} t }\sum\limits_{a=1}^{\infty} \frac{1}{\cosh(\frac{\pi da^2}{2t})} &\\&+ \frac{t\sinh(\pi t)G_{-N,t,k}}{2\sqrt{k}\pi}\sum\limits_{r=0}^{\infty} \sum\limits_{a=1}^{\infty}\frac{(-1)^{r-1}(2r+1)e^{-\pi t(2r+1)}}{t^2(2r+1)^2+d^2a^4} +\frac{\sinh(\pi t)}{8\sqrt{k}t}\sum\limits_{r=1}^{\infty}\sum\limits_{a=1}^{\infty}\frac{(-1)^rG_{r-N,t,k}}{ \cosh(\frac{\pi (r-da^2)}{2t})} &\\&+\frac{\sinh(\pi t)}{8\sqrt{k}t}\sum\limits_{r=1}^{\infty}\sum\limits_{a=1}^{\infty}\frac{(-1)^rG_{-r-N,t,k}}{ \cosh(\frac{\pi (r+da^2)}{2t})}  -\frac{t\sinh(\pi t)}{4\sqrt{k} \pi}\sum\limits_{r=1}^{\infty} \sum\limits_{m=0}^{\infty}   \frac{(-1)^{m-1}(2m+1)e^{-\pi t(2m+1)}G_{r-N,t,k}}{t^2(2m+1)^2+r^2} &\\&-\frac{t\sinh(\pi t)}{4\sqrt{k}\pi}\sum\limits_{r=1}^{\infty} \sum\limits_{m=0}^{\infty} \frac{(-1)^{m-1}(2m+1)e^{-\pi t(2m+1)}G_{-r-N,t,k}}{t^2(2m+1)^2+r^2} &\\&+\frac{t\pi\sinh(\pi t)}{8d^2\sqrt{k} \pi}\sum\limits_{r=1}^{\infty} \sum\limits_{m=0}^{\infty} (-1)^{m-1}(2m+1)e^{-\pi t(2m+1)}G_{r-N,t,k} T_{-\frac{r}{d},\frac{t(2m+1)}{d}} &\\& +\frac{\pi t\sinh(\pi t)}{8d^2\sqrt{k}\pi}\sum\limits_{r=1}^{\infty} \sum\limits_{m=0}^{\infty}(-1)^{m-1}(2m+1)e^{-\pi t(2m+1)}G_{-r-N,t,k}T_{\frac{r}{d}, \frac{t(2m+1)}{d}}  ,\label{evenFormulaSumofSquares} \tag{1.9}
\end{flalign*}
when $d$ is even, and

\begin{flalign*}
    &\frac{1}{k^2}\sum\limits_{\substack{da^2+kb^2=N \\ (a,b) \in \mathbb{N}^2}} \frac{1}{b^4} 
    \\ &= \left(\frac{ \pi^2}{12} - \frac{\pi \coth(\pi t)}{4} + \frac{\pi \coth(\pi t)}{48 t^2} \right)\epsilon_N + \frac{\pi^2}{3k(e^{2\pi t}-1)}{\sum\limits_{a=1}^{\infty} }^*\frac{1}{N-da^2}+ \frac{\pi - \coth(\pi t)}{2\pi}{\sum\limits_{a=1}^{\infty} }^*\frac{1}{(N-da^2)^2} &\\&+ \frac{(-1)^{N}(1-t^2)\coth(\pi t)}{\pi}{\sum\limits_{a=1}^{\infty} }^* \frac{(-1)^{a}}{(N-da^2)^2} + \frac{(-1)^{N} \pi^3 \coth(\pi t)}{12kt} {\sum\limits_{a=1}^{\infty}}^* \frac{(-1)^{a-1}}{\sinh(\frac{\pi (N-da^2)}{2t})}  &\\&+ \frac{(-1)^{N}\coth(\pi t)(t-\frac{1}{t})}{4} {\sum\limits_{a=1}^{\infty}}^* \frac{(-1)^{a}}{(N-da^2)\sinh(\frac{\pi (N-da^2)}{2t})} + \frac{(-1)^{N}\pi \coth(\pi t)}{8} {\sum\limits_{a=1}^{\infty}}^* \frac{(-1)^{a}\cosh(\frac{\pi(N-da^2)}{2t})}{\sinh^2(\frac{\pi (N-da^2)}{2t})} &\end{flalign*} \begin{flalign*}&- \frac{\pi^2\coth(\pi t)}{3k} \sum\limits_{r=1}^{\infty} \sum\limits_{a=1}^{\infty} \frac{(-1)^{r-1}e^{-2\pi t r}(N-da^2)}{4t^2r^2+(N-da^2)^2}  - t\coth(\pi t) \sum\limits_{r=1}^{\infty} \sum\limits_{a=1}^{\infty} \frac{(-1)^{r-1}re^{-2\pi t r}}{4t^2r^2+(N-da^2)^2} &\\& - \frac{\coth(\pi t)}{\pi}\sum\limits_{r=1}^{\infty}\sum\limits_{a=1}^{\infty} \frac{(-1)^{r-1}e^{-2\pi t r}(4t^2r^2-(N-da^2)^2)}{(4t^2r^2+(N-da^2)^2)^2} + \frac{\sinh(\pi t) G_{-N,t,k}}{8\sqrt{k} t }\sum\limits_{a=1}^{\infty} \frac{1}{\cosh(\frac{\pi da^2}{2t})} &\\&+ \frac{t\sinh(\pi t)G_{-N,t,k}}{2\sqrt{k}\pi}\sum\limits_{r=0}^{\infty} \sum\limits_{a=1}^{\infty}\frac{(-1)^{r+a-1}(2r+1)e^{-\pi t(2r+1)}}{t^2(2r+1)^2+d^2a^4} +\frac{\sinh(\pi t)}{8\sqrt{k}t}\sum\limits_{r=1}^{\infty}\sum\limits_{a=1}^{\infty}\frac{(-1)^rG_{r-N,t,k}}{ \cosh(\frac{\pi (r-da^2)}{2t})}&\\&+\frac{\sinh(\pi t)}{8\sqrt{k}t}\sum\limits_{r=1}^{\infty}\sum\limits_{a=1}^{\infty}\frac{(-1)^rG_{-r-N,t,k}}{ \cosh(\frac{\pi (r+da^2)}{2t})} -\frac{t\sinh(\pi t)}{4\sqrt{k} \pi}\sum\limits_{r=1}^{\infty} \sum\limits_{m=0}^{\infty}   \frac{(-1)^{m-1}(2m+1)e^{-\pi t(2m+1)}G_{r-N,t,k}}{t^2(2m+1)^2+r^2}&\\&-\frac{t\sinh(\pi t)}{4\sqrt{k}\pi}\sum\limits_{r=1}^{\infty} \sum\limits_{m=0}^{\infty} \frac{(-1)^{m-1}(2m+1)e^{-\pi t(2m+1)}G_{-r-N,t,k}}{t^2(2m+1)^2+r^2} &\\& +\frac{\pi t\sinh(\pi t)}{4d^2\sqrt{k} \pi}\sum\limits_{r=1}^{\infty} \sum\limits_{m=0}^{\infty} (-1)^{m-1}(2m+1)e^{-\pi t(2m+1)}G_{r-N,t,k} V_{-\frac{r}{d},\frac{t(2m+1)}{d}} &\\&+\frac{\pi t\sinh(\pi t)}{4d^2\sqrt{k}\pi}\sum\limits_{r=1}^{\infty} \sum\limits_{m=0}^{\infty}(-1)^{m-1}(2m+1)e^{-\pi t(2m+1)}G_{-r-N,t,k}  V_{\frac{r}{d}, \frac{t(2m+1)}{d}} ,\label{oddFormulaSumofSquares} \tag{1.10}
\end{flalign*}
when $d$ is odd. Where
\[
\epsilon_N = \begin{cases}
    1&, \textit{ if } N = dm^2 \textit{ for some natural number }m, \\
    0&, \textit{ if else,}
\end{cases}
\]
and the $*$ indicates that the term inside the summation exhibiting a singularity is taken as zero.
\end{proposition}

Next, we state the counterpart for the sum whose terms are summed over the positive integer solutions $(a,b)$ of the Diophantine equation $kb^2-da^2 = N$.
\vspace{0.7cm}
\begin{proposition} For natural numbers $N, d$, and $k$, and for real value $t>0$, we have 
\begin{flalign*}
    &\frac{1}{k^2}\sum\limits_{\substack{kb^2-da^2=N \\ (a,b) \in \mathbb{N}^2}} \frac{1}{b^4}  &\\& =  \frac{\pi^2}{3k(e^{2\pi t}-1)} \sum\limits_{a=1}^{\infty} \frac{1}{N+da^2}+ \frac{\pi - \coth(\pi t)}{2\pi}\sum\limits_{a=1}^{\infty} \frac{1}{(N+da^2)^2}+ \frac{(-1)^{N}(1-t^2)\coth(\pi t)}{\pi }\sum\limits_{a=1}^{\infty} \frac{1}{(N+da^2)^2} &\\& + \frac{(-1)^{N-1} \pi^3 \coth(\pi t)}{12kt }\sum\limits_{a=1}^{\infty} \frac{1}{\sinh(\frac{\pi (N+da^2)}{2t})}  + \frac{(-1)^{N}\coth(\pi t)(t-\frac{1}{t})}{4}\sum\limits_{a=1}^{\infty} \frac{1}{(N+da^2)\sinh(\frac{\pi (N+da^2)}{2t})}&\\&+ \frac{(-1)^{N}\pi \coth(\pi t)}{8}\sum\limits_{a=1}^{\infty} \frac{\cosh(\frac{\pi(N+da^2)}{2t})}{\sinh^2(\frac{\pi (N+da^2)}{2t})} - \frac{\pi^2\coth(\pi t)}{3k} \sum\limits_{r=1}^{\infty}\sum\limits_{a=1}^{\infty} \frac{(-1)^{r-1}e^{-2\pi t r}(N+da^2)}{4t^2r^2+(N+da^2)^2} &\\& - t\coth(\pi t) \sum\limits_{r=1}^{\infty} \sum\limits_{a=1}^{\infty} \frac{(-1)^{r-1}re^{-2\pi t r}}{4t^2r^2+(N+da^2)^2} - \frac{\coth(\pi t)}{\pi}\sum\limits_{r=1}^{\infty} \sum\limits_{a=1}^{\infty} \frac{(-1)^{r-1}e^{-2\pi t r}(4t^2r^2-(N+da^2)^2)}{(4t^2r^2+(N+da^2)^2)^2} &\\& + \frac{\sinh(\pi t) G_{-N,t,k}}{8\sqrt{k} t } \sum\limits_{a=1}^{\infty} \frac{1}{\cosh(\frac{\pi da^2}{2t})} + \frac{t\sinh(\pi t)G_{-N,t,k}}{2\sqrt{k}\pi}\sum\limits_{r=0}^{\infty} \sum\limits_{a=1}^{\infty}\frac{(-1)^{r-1}(2r+1)e^{-\pi t(2r+1)}}{t^2(2r+1)^2+d^2a^4} &\\&+\frac{\sinh(\pi t)}{8\sqrt{k}t}\sum\limits_{r=1}^{\infty} \sum\limits_{a=1}^{\infty}\frac{(-1)^rG_{r-N,t,k}}{ \cosh(\frac{\pi (r+da^2)}{2t})} +\frac{\sinh(\pi t)}{8\sqrt{k}t}\sum\limits_{r=1}^{\infty} \sum\limits_{a=1}^{\infty}(-1)^r\frac{G_{-r-N,t,k}}{ \cosh(\frac{\pi (r-da^2)}{2t})} &\\&-\frac{t\sinh(\pi t)}{4\sqrt{k} \pi}\sum\limits_{r=1}^{\infty} \sum\limits_{m=0}^{\infty}\frac{(-1)^{m-1}(2m+1)e^{-\pi t(2m+1)}G_{r-N,t,k}}{t^2(2m+1)^2+ r^2} -\frac{t\sinh(\pi t)}{4\sqrt{k}\pi}\sum\limits_{r=1}^{\infty} \sum\limits_{m=0}^{\infty} \frac{(-1)^{m-1}(2m+1)e^{-\pi t(2m+1)}G_{-r-N,t,k}}{t^2(2m+1)^2+ r^2}&\\& +\frac{\pi t\sinh(\pi t)}{8d^2\sqrt{k} \pi}\sum\limits_{r=1}^{\infty} \sum\limits_{m=0}^{\infty}(-1)^{m-1}(2m+1)e^{-\pi t(2m+1)}G_{r-N,t,k} T_{\frac{r}{d}, \frac{t(2m+1)}{d}} &\\&+\frac{\pi t\sinh(\pi t)}{8d^2\sqrt{k}\pi}\sum\limits_{r=1}^{\infty} \sum\limits_{m=0}^{\infty}(-1)^{m-1}(2m+1)e^{-\pi t(2m+1)}G_{-r-N,t,k} T_{-\frac{r}{d}, \frac{t(2m+1)}{d}}  ,\label{evenFormulaDifferenceofSquares} \tag{1.11}
\end{flalign*}
when $d$ is even, and

\begin{flalign*}
    &\frac{1}{k^2}\sum\limits_{\substack{kb^2-da^2=N \\ (a,b) \in \mathbb{N}^2}} \frac{1}{b^4}  &\\& =  \frac{\pi^2}{3k(e^{2\pi t}-1)} \sum\limits_{a=1}^{\infty} \frac{1}{N+da^2}+ \frac{\pi - \coth(\pi t)}{2\pi}\sum\limits_{a=1}^{\infty} \frac{1}{(N+da^2)^2}+ \frac{(-1)^{N}(1-t^2)\coth(\pi t)}{\pi }\sum\limits_{a=1}^{\infty} \frac{(-1)^{a}}{(N+da^2)^2} &\\& + \frac{(-1)^{N} \pi^3 \coth(\pi t)}{12kt }\sum\limits_{a=1}^{\infty} \frac{(-1)^{a-1}}{\sinh(\frac{\pi (N+da^2)}{2t})}  + \frac{(-1)^{N}\coth(\pi t)(t-\frac{1}{t})}{4}\sum\limits_{a=1}^{\infty} \frac{(-1)^{a}}{(N+da^2)\sinh(\frac{\pi (N+da^2)}{2t})}&\\&+ \frac{(-1)^{N}\pi \coth(\pi t)}{8}\sum\limits_{a=1}^{\infty} \frac{(-1)^{a}\cosh(\frac{\pi(N+da^2)}{2t})}{\sinh^2(\frac{\pi (N+da^2)}{2t})} - \frac{\pi^2\coth(\pi t)}{3k} \sum\limits_{r=1}^{\infty}\sum\limits_{a=1}^{\infty} \frac{(-1)^{r-1}e^{-2\pi t r}(N+da^2)}{4t^2r^2+(N+da^2)^2} &\\& - t\coth(\pi t) \sum\limits_{r=1}^{\infty} \sum\limits_{a=1}^{\infty} \frac{(-1)^{r-1}re^{-2\pi t r}}{4t^2r^2+(N+da^2)^2} - \frac{\coth(\pi t)}{\pi}\sum\limits_{r=1}^{\infty} \sum\limits_{a=1}^{\infty} \frac{(-1)^{r-1}e^{-2\pi t r}(4t^2r^2-(N+da^2)^2)}{(4t^2r^2+(N+da^2)^2)^2} &\\& + \frac{\sinh(\pi t) G_{-N,t,k}}{8\sqrt{k} t } \sum\limits_{a=1}^{\infty} \frac{1}{\cosh(\frac{\pi da^2}{2t})} + \frac{t\sinh(\pi t)G_{-N,t,k}}{2\sqrt{k}\pi}\sum\limits_{r=0}^{\infty} \sum\limits_{a=1}^{\infty}\frac{(-1)^{r+a-1}(2r+1)e^{-\pi t(2r+1)}}{t^2(2r+1)^2+d^2a^4} &\\&+\frac{\sinh(\pi t)}{8\sqrt{k}t}\sum\limits_{r=1}^{\infty} \sum\limits_{a=1}^{\infty}\frac{(-1)^rG_{r-N,t,k}}{ \cosh(\frac{\pi (r+da^2)}{2t})} +\frac{\sinh(\pi t)}{8\sqrt{k}t}\sum\limits_{r=1}^{\infty} \sum\limits_{a=1}^{\infty}(-1)^r\frac{G_{-r-N,t,k}}{ \cosh(\frac{\pi (r-da^2)}{2t})} &\\&-\frac{t\sinh(\pi t)}{4\sqrt{k} \pi}\sum\limits_{r=1}^{\infty} \sum\limits_{m=0}^{\infty} \frac{(-1)^{m-1}(2m+1)e^{-\pi t(2m+1)}G_{r-N,t,k}}{t^2(2m+1)^2+r^2} -\frac{t\sinh(\pi t)}{4\sqrt{k}\pi}\sum\limits_{r=1}^{\infty} \sum\limits_{m=0}^{\infty} \frac{(-1)^{m-1}(2m+1)e^{-\pi t(2m+1)}G_{-r-N,t,k}}{t^2(2m+1)^2+r^2}
    &\\&+\frac{\pi t\sinh(\pi t)}{4d^2\sqrt{k} \pi}\sum\limits_{r=1}^{\infty} \sum\limits_{m=0}^{\infty}(-1)^{m-1}(2m+1)e^{-\pi t(2m+1)}G_{r-N,t,k} V_{\frac{r}{d}, \frac{t(2m+1)}{d}} &\\&+\frac{\pi t\sinh(\pi t)}{4d^2\sqrt{k}\pi}\sum\limits_{r=1}^{\infty} \sum\limits_{m=0}^{\infty}(-1)^{m-1}(2m+1)e^{-\pi t(2m+1)}G_{-r-N,t,k}V_{-\frac{r}{d}, \frac{t(2m+1)}{d}}, \label{oddFormulaDifferenceofSquares} \tag{1.12}
\end{flalign*}
when $d$ is odd.
    
\end{proposition}

As a first step towards our demonstration, we shall perform an inversion of a general series, provided in the next section.
\vspace{0.7cm}
\section{Inversion of a general series}

Now, we go on to describe a certain kind of invertibility of the following formal series of partial fractions  

\[
\sum\limits_{n=1}^{\infty} \frac{f(n)}{n+z}. \label{2parfrac-series} \tag{2.1}
\]

The following lemma will be necessary.
\vspace{0.7cm}
\begin{lemma}
    If for a non-negative sequence of real numbers $\{f(n)\}_{n\in \mathbb{N}}$ and arbitrary complex value $z$, we have
    \[ \sum\limits_{n=1}^{\infty} \frac{f(n)}{n+z} = F(z).\label{parfrac-equality} \tag{2.2}\]
    Then for $|\beta|\leq 1$, and real $t$ there holds
    \begin{multline} 
        \frac{\pi \cosh(\pi\beta t)}{\sinh(\pi t)} \sum\limits_{n=1}^{\infty} (-1)^n e^{-\pi i n\beta} f(n)  \\=\frac{F(-it)-F(it)}{2i}+\frac{1}{2i}\sum\limits_{k=1}^{\infty} (-1)^ke^{\pi i k \beta}(F(k-it)-F(k+it))+\frac{1}{2i}\sum\limits_{k=1}^{\infty} (-1)^ke^{-\pi i k \beta}(F(-k-it)-F(-k+it))  . 
        \label{parfrac-coefficient-series} \tag{2.3}
    \end{multline} 
    Provided that all the series from the right hand side converge absolutely.
\end{lemma}

\begin{proof}[Proof of Lemma 4] Suppose that \eqref{parfrac-equality} holds, then we have
    \[
    \sum\limits_{n=1}^{\infty} \frac{f(n)}{(n+z)^2+t^2}=\frac{F(z-it)-F(z+it)}{2it}. \label{proofeq1} \tag{2.4}
    \]

    Therefore, we see that

    \begin{multline}
        \frac{F(-it)-F(it)}{2it}+\frac{1}{2it}\sum\limits_{k=1}^{\infty} (-1)^ke^{\pi i k \beta}(F(k-it)-F(k+it))+\frac{1}{2it}\sum\limits_{k=1}^{\infty} (-1)^ke^{-\pi i k \beta}(F(-k-it)-F(-k+it)) \\ 
        = \sum\limits_{k=0}^{\infty} \sum\limits_{n=1}^{\infty} \frac{(-1)^ke^{\pi i k\beta}f(n)}{(n+k)^2+t^2}+\sum\limits_{k=1}^{\infty} \sum\limits_{n=1}^{\infty} \frac{(-1)^ke^{-\pi i k\beta}f(n)}{(n-k)^2+t^2}. \label{proofeq2} \tag{2.5}
    \end{multline}

    Where it is assumed that the series $ \sum\limits_{k=1}^{\infty} (-1)^ke^{\pi i k \beta}(F(k-it)-F(k+it))$ and $\sum\limits_{k=1}^{\infty} (-1)^ke^{-\pi i k \beta}(F(-k-it)$\\$-F(-k+it)) $ being absolutely convergent. Now in each double sum in the right hand side, the order of summation can be interchanged for all real $t$ since each of the series converges absolutely as series over $k$ and each of the terms $\frac{f(n)}{(n+k)^2+t^2}$ and $\frac{f(n)}{(n-k)^2+t^2}$ is non-negative for real $t$. Thus, we have

    \begin{flalign*}
        &\frac{F(-it)-F(it)}{2it}+\frac{1}{2it}\sum\limits_{k=1}^{\infty} (-1)^ke^{\pi i k \beta}(F(k-it)-F(k+it))+\frac{1}{2it}\sum\limits_{k=1}^{\infty} (-1)^ke^{-\pi i k \beta}(F(-k-it)-F(-k+it)) & \\ 
        &= \sum\limits_{n=1}^{\infty} \sum\limits_{k=0}^{\infty} \frac{(-1)^ke^{\pi i k\beta}f(n)}{(n+k)^2+t^2}+\sum\limits_{n=1}^{\infty} \sum\limits_{k=1}^{\infty} \frac{(-1)^ke^{-\pi i k\beta}f(n)}{(n-k)^2+t^2} 
        &\\ &= \sum\limits_{n=1}^{\infty} (-1)^n e^{-\pi i n\beta} f(n)\sum\limits_{k=0}^{\infty} \frac{(-1)^{n+k}e^{\pi i (n+k)\beta}}{(n+k)^2+t^2}+\sum\limits_{n=1}^{\infty}(-1)^ne^{-\pi i n\beta}f(n) \sum\limits_{k=1}^{\infty} \frac{(-1)^{n-k}e^{\pi i (n-k)\beta}}{(n-k)^2+t^2} 
        &\\& = \sum\limits_{n=1}^{\infty} (-1)^n e^{-\pi i n\beta} f(n)\sum\limits_{k=0}^{\infty} \frac{(-1)^{n+k}e^{\pi i (n+k)\beta}}{(n+k)^2+t^2}+\sum\limits_{n=1}^{\infty}(-1)^ne^{-\pi i n\beta}f(n) \sum\limits_{k=1}^{\infty} \frac{(-1)^{n-k}e^{\pi i (n-k)\beta}}{(n-k)^2+t^2} &\\& 
        = \sum\limits_{n=1}^{\infty} (-1)^n e^{-\pi i n\beta} f(n) \left\{ \sum\limits_{k=0}^{\infty} \frac{(-1)^{n+k}e^{\pi i (n+k)\beta}}{(n+k)^2+t^2} + \sum\limits_{k=1}^{\infty} \frac{(-1)^{n-k}e^{\pi i (n-k)\beta}}{(n-k)^2+t^2}\right\} & \\& 
        = \sum\limits_{n=1}^{\infty} (-1)^n e^{-\pi i n\beta} f(n) \left\{ \sum\limits_{k=0}^{\infty} \frac{(-1)^{k}e^{\pi i k\beta}}{k^2+t^2} + \sum\limits_{k=1}^{\infty} \frac{(-1)^{k}e^{-\pi i k\beta}}{k^2+t^2}\right\} & 
        \label{proofeq3} \tag{2.6}
    \end{flalign*}

Now we finally note that for $|\beta| \leq 1$,  $\sum\limits_{k=0}^{\infty}\frac{(-1)^{k}e^{\pi i k\beta}}{k^2+t^2} + \sum\limits_{k=1}^{\infty} \frac{(-1)^{k}e^{-\pi i k\beta}}{k^2+t^2} = \frac{\pi \cosh(\pi \beta t)}{t \sinh(\pi t)}$, and the proof is complete.
    
\end{proof}

The following theorem performs a kind of inversion on the series \eqref{2parfrac-series}. 
\vspace{0.7cm}
\begin{theorem}
    If for a non-negative sequence of real numbers $\{f(n)\}_{n\in \mathbb{N}}$ and arbitrary complex value $z$, we have
    \[ \sum\limits_{n=1}^{\infty} \frac{f(n)}{n+z} = F(z). \label{2parfrac-equality} \tag{2.7}\]
    Then for each natural number $N$, and real $t \neq 0$

    \begin{flalign*}
        &   f(N)  &\\& =  \frac{\sinh(\pi t)}{2i\pi }\sum\limits_{k=1}^{\infty} (-1)^{k}(F(k-N-it)-F(k-N+it)+F(-k-N-it)-F(-k-N+it))\int_{0}^{1} \frac{\cos(\pi k \beta)}{\cosh(\pi\beta t)} d\beta &\\& +\frac{(F(-N-it)-F(-N+it))\sinh(\pi t)\arctan(\tanh(\frac{\pi t}{2}))}{i\pi^2t}  .\label{parfrac-inverted-equality} \tag{2.8}
    \end{flalign*}
    Provided that the series $\sum\limits_{k=1}^{\infty} (-1)^ke^{\pi i k \beta}(F(k-it)-F(k+it))$ and $\sum\limits_{k=1}^{\infty} (-1)^ke^{-\pi i k \beta}(F(-k-it)-F(-k+it))$ converges absolutely and uniformly for $|\beta| \leq 1$. The expression in right hand side of \eqref{parfrac-inverted-equality} vanishes when $N$ is either zero or a negative integer.
    
\end{theorem}

\begin{proof}[Proof of Theorem 5]
Since \eqref{2parfrac-equality} holds, by \textbf{Lemma 4} we have
    \begin{flalign*}
        &\frac{\pi \cosh(\pi\beta t)}{\sinh(\pi t)} \sum\limits_{n=1}^{\infty} (-1)^n e^{-\pi i n\beta} f(n)  \\&=\frac{F(-it)-F(it)}{2i}+\frac{1}{2i}\sum\limits_{k=1}^{\infty} (-1)^ke^{\pi i k \beta}(F(k-it)-F(k+it))+\frac{1}{2i}\sum\limits_{k=1}^{\infty} (-1)^ke^{-\pi i k \beta}(F(-k-it)-F(-k+it))  .  
    \end{flalign*}

Where the absolute convergence of the series from the right hand side is assumed. When we separate the real and imaginary parts from both sides, we get the equalities

\begin{flalign*}
        & \sum\limits_{n=1}^{\infty} (-1)^n \sin(\pi n\beta) f(n) &\\& =   \frac{i\sinh(\pi t)}{2\pi \cosh(\pi\beta t)}\sum\limits_{k=1}^{\infty} (-1)^k\sin(\pi  k \beta)(F(k-it)-F(k+it)-F(-k-it)+F(-k+it))
        \label{proofeq4} \tag{2.9}
\end{flalign*}

\begin{flalign*}
        & \sum\limits_{n=1}^{\infty} (-1)^n \cos(\pi n\beta) f(n)   &\\&=\frac{(F(-it)-F(it))\sinh(\pi t)}{2i\pi\cosh(\pi\beta t)}+\frac{\sinh(\pi t)}{2i\pi \cosh(\pi \beta t)}\sum\limits_{k=1}^{\infty} (-1)^k \cos(\pi  k \beta)(F(k-it)-F(k+it)+F(-k-it)-F(-k+it))  . \label{proofeq5} \tag{2.10} 
    \end{flalign*}

Thus in each equality, we can recover $(-1)^Nf(N)$ for each natural number $N$, as the Fourier coefficients \cite[p.~363]{ComplexFunctionsReinhold} of the expression from the right hand side

\begin{flalign*}
        & (-1)^N f(N) &\\& =   \int_{-1}^{1}\frac{i\sin(\pi N \beta)\sinh(\pi t)}{2\pi \cosh(\pi\beta t)}\left\{\sum\limits_{k=1}^{\infty} (-1)^k\sin(\pi  k \beta)(F(k-it)-F(k+it)-F(-k-it)+F(-k+it)) \right\} d\beta &\\& =\int_{-1}^{1}\left\{\frac{(F(-it)-F(it))\cos(\pi N \beta)\sinh(\pi t)}{2i\pi\cosh(\pi\beta t)} \right. &\\& \left. +\frac{\cos(\pi N \beta)\sinh(\pi t)}{2i\pi \cosh(\pi \beta t)}\sum\limits_{k=1}^{\infty} (-1)^k \cos(\pi  k \beta)(F(k-it)-F(k+it)+F(-k-it)-F(-k+it))  \right\} d\beta.& \label{proofeq6} \tag{2.11}  \end{flalign*}

When the uniform convergence of the series $\sum\limits_{k=1}^{\infty} (-1)^ke^{\pi i k \beta}(F(k-it)-F(k+it))$ and $\sum\limits_{k=1}^{\infty} (-1)^ke^{-\pi i k \beta}(F(-k-it)-F(-k+it))$ for $|\beta| \leq 1$ is assumed, the expressions in \eqref{proofeq6} can be integrated term-wise, thus yielding the equalities. 

\begin{flalign*}
        &   f(N) &\\& =   \frac{i\sinh(\pi t)}{\pi }\sum\limits_{k=1}^{\infty} (-1)^{N+k}(F(k-it)-F(k+it)-F(-k-it)+F(-k+it)) \int_{0}^{1} \frac{\sin(\pi N \beta)\sin(\pi  k \beta)}{\cosh(\pi\beta t)} d\beta .\label{proofeq7} \tag{2.12}
    \end{flalign*}

\begin{flalign*}
        &   f(N) &\\& = -\frac{i\sinh(\pi t)}{\pi}\sum\limits_{k=1}^{\infty} (-1)^{N+k} (F(k-it)-F(k+it)+F(-k-it)-F(-k+it))\int_{0}^{1} \frac{\cos(\pi N \beta)\cos(\pi  k \beta)}{\cosh(\pi \beta t)} d\beta &\\& +(-1)^N\frac{(F(-it)-F(it))\sinh(\pi t)}{i\pi} \int_{0}^{1}\frac{\cos(\pi N \beta)}{ \cosh(\pi\beta t)} d\beta .\label{proofeq8} \tag{2.13}
    \end{flalign*}

Addition of the two equalities gives

\begin{flalign*}
        &   2f(N) &\\& =  \frac{\sinh(\pi t)}{i\pi }\sum\limits_{k=1}^{\infty} (-1)^{N+k}(F(k-it)-F(k+it))\int_{0}^{1} \frac{\cos(\pi (k+N) \beta)}{\cosh(\pi\beta t)} d\beta &\\& +\frac{\sinh(\pi t)}{i\pi}\sum\limits_{k=1}^{\infty}(-1)^{N+k}(F(-k-it)-F(-k+it))\int_{0}^{1}\frac{\cos(\pi (-k+N) \beta)}{ \cosh(\pi\beta t)} d\beta &\\&+(-1)^N\frac{(F(-it)-F(it))\sinh(\pi t)}{i\pi} \int_{0}^{1}\frac{\cos(\pi N \beta)}{ \cosh(\pi\beta t)} d\beta &\end{flalign*} \begin{flalign*}& =  \frac{\sinh(\pi t)}{i\pi }\sum\limits_{k=1}^{\infty} (-1)^{k}(F(k-N-it)-F(k-N+it))\int_{0}^{1} \frac{\cos(\pi k \beta)}{\cosh(\pi\beta t)} d\beta &\\& +\frac{\sinh(\pi t)}{i\pi}\sum\limits_{k=1}^{\infty}(-1)^{k}(F(-k-N-it)-F(-k-N+it))\int_{0}^{1}\frac{\cos(\pi k \beta)}{ \cosh(\pi\beta t)} d\beta &\\&+\frac{(F(-N-it)-F(-N+it))\sinh(\pi t)}{i\pi} \int_{0}^{1}\frac{1}{ \cosh(\pi\beta t)} d\beta .\label{proofeq9} \tag{2.14}
    \end{flalign*}

    Noting that $\int_{0}^{1} \sech(\pi\beta t) d\beta = \frac{2\arctan(\tanh(\frac{\pi t}{2}))}{\pi t}$, we arrive at the desired equality. The fact that this last expression vanishes if $N$ is either $0$ or a negative integer follows from the observation that when $N=0$ both of the expressions in the right hand side of \eqref{proofeq7} and \eqref{proofeq8} are $0$, and the case for negative integers $N$ follows from the subtraction of the two equalities \eqref{proofeq7} and \eqref{proofeq8}.
    
\end{proof}

Then it readily follows from the case $N=0$ of \textbf{Theorem 5} that if $F(z)$ is a function which can be given an expansion of the form in \eqref{2parfrac-equality}, and if it happens to hold that the series $\sum\limits_{k=1}^{\infty} (-1)^ke^{\pi i k \beta}(F(k-it)-F(k+it))$ and $\sum\limits_{k=1}^{\infty} (-1)^ke^{-\pi i k \beta}(F(-k-it)-F(-k+it))$ converges absolutely and uniformly for $|\beta| \leq 1$, we will always have the equality 

\begin{flalign*}
        &   \frac{\arctan(\tanh(\frac{\pi t}{2}))}{\pi t}(F(it)-F(-it))  &\\& =  \frac{1}{2 }\sum\limits_{k=1}^{\infty} (-1)^{k}(F(k-it)-F(k+it)+F(-k-it)-F(-k+it))\int_{0}^{1} \frac{\cos(\pi k \beta)}{\cosh(\pi\beta t)} d\beta . \label{seriesIdentityForF} \tag{2.15}
    \end{flalign*}
\vspace{0.7cm}
\section{The arithmetic function $q_{k,s}(N)$}

Here we propose to consider a certain class of arithmetic functions, defined for a pair of natural numbers $k$ and $s$. Then for each such a pair, we define the arithmetic function $q_{k,s}(N)$ for all natural numbers $N$ in the following manner

\[
q_{k,s}(N) = \begin{cases}
    1 & \textit{ if } N = km^{2s} \textit{ for some natural number } m, \\
    0 & \textit{ if else}.
         \end{cases} \label{qDef} \tag{3.1}
\]

In our current section, we shall demonstrate that $q_{k,s}(N)$ can always be given an analytic expression in the form of a convergent infinite series. We shall first consider the case for $s=1$. We shall denote $q_{k,1}(N)$ equivalently as $q_k(N)$. Then we remark that $q_k(N)=1$ only when $N$ is divisible by $k$ and $\frac{N}{k}$ is a square, and vanishes in all other cases. Firstly we shall state the following lemma.
\vspace{0.7cm}
\begin{lemma} For arbitrary complex value $z$, we have
    \[
    \sum\limits_{n=1}^{\infty} \frac{q_k(n)}{n^2(n+z)} = \frac{\pi^4}{90k^2z}- \frac{\pi^2}{6kz^2} - \frac{1}{2z^3} + \frac{\pi\cosh(\frac{\pi \sqrt{z}}{\sqrt{k}})}{2z^2\sqrt{zk}\sinh(\frac{\pi \sqrt{z}}{\sqrt{k}})}. \label{generatingSeriesForq} \tag{3.2}
    \]
\end{lemma}

\begin{proof}[Proof of Lemma 6]
    First we make the observation that 
    \begin{flalign*}
        \sum\limits_{n=1}^{\infty} \frac{q_k(n)}{n^2(n+z)} &= \sum\limits_{n=1}^{\infty} \frac{1}{k^2n^4(kn^2+z)}& \\
        &=  \frac{1}{k^2z}\sum\limits_{n=1}^{\infty} \frac{1}{n^4}- \frac{1}{kz^2} \sum\limits_{n=1}^{\infty} \frac{1}{n^2} + \frac{1}{kz^2}\sum\limits_{n=1}^{\infty} \frac{1}{n^2+\frac{z}{k}}&  \label{proofeq10} \tag{3.3}
    \end{flalign*}

    And the desired equality follows from the known values of $\zeta(4)$, and $\zeta(2)$, and the Mittag-Leffler expansion of the hyperbolic cotangent function.
\end{proof}

From \textbf{Lemma 6} we see that if $F(z) = \frac{\pi^4}{90k^2z}- \frac{\pi^2}{6kz^2} - \frac{1}{2z^3} + \frac{\pi\cosh(\frac{\pi \sqrt{z}}{\sqrt{k}})}{2z^2\sqrt{zk}\sinh(\frac{\pi \sqrt{z}}{\sqrt{k}})}$, satisfies the requirements described by \textbf{Theorem 5}, we can obtain an analytic representation of $\frac{q_k(N)}{N^2}$ in the form of a convergent infinite series. Thus we shall try to prove that for $F(z)$ defined as above, the series $\sum\limits_{r=1}^{\infty} (-1)^re^{\pi i r \beta}(F(r-it)-F(r+it))$ and $\sum\limits_{r=1}^{\infty} (-1)^re^{-\pi i r \beta}(F(-r-it)-F(-r+it))$ converges absolutely and uniformly for $|\beta| \leq 1$.
\vspace{0.7cm}
\begin{lemma} If $F(z)$ is given by
\[F(z)= \frac{\pi^4}{90k^2z}- \frac{\pi^2}{6kz^2} - \frac{1}{2z^3} + \frac{\pi\cosh(\frac{\pi \sqrt{z}}{\sqrt{k}})}{2z^2\sqrt{zk}\sinh(\frac{\pi \sqrt{z}}{\sqrt{k}})},\label{2generatingSeriesForq} \tag{3.4}\]

and $H(z)= F(z-it)-F(z+it)$, then the series $\sum\limits_{r=1}^{\infty} (-1)^re^{\pi i r \beta}H(r)$ and $\sum\limits_{k=1}^{\infty} (-1)^re^{-\pi i r \beta}H(-r)$ converges absolutely and uniformly for all natural numbers $k$, and $|\beta| \leq 1$.
    
\end{lemma}
\begin{proof}[Proof of Lemma 7]
   Let us first compute $H(M) = F(M-it)-F(M+it)$ 
   \begin{flalign*}
       F(M-it)-F(M+it) = & \frac{\pi^4}{90k^2(M-it)}-\frac{\pi^4}{90k^2(M+it)}- \frac{\pi^2}{6k(M-it)^2}+\frac{\pi^2}{6k(M+it)^2} &\\&- \frac{1}{2(M-it)^3}+\frac{1}{2(M+it)^3} + \frac{\pi\cosh(\frac{\pi \sqrt{M-it}}{\sqrt{k}})}{2(M-it)^2\sqrt{(M-it)k}\sinh(\frac{\pi \sqrt{M-it}}{\sqrt{k}})} &\\& - \frac{\pi\cosh(\frac{\pi \sqrt{M+it}}{\sqrt{k}})}{2(M+it)^2\sqrt{(M+it)k}\sinh(\frac{\pi \sqrt{M+it}}{\sqrt{k}})}.& \label{proofeq11} \tag{3.5}
   \end{flalign*}
If $u_{M,t}$ and $v_{M,t}$ are defined by \eqref{uvdef}, we have

\begin{flalign*}
       &H(M)  &\\& = \frac{it\pi^4}{30k^2(M^2+t^2)}- \frac{2\pi^2itM}{3k(M^2+t^2)^2} - it\frac{(3M^2-t^2)}{(M^2+t^2)^3} &\\& + \frac{i\pi((M^2u_{M,t}-t^2u_{M,t}-2Mtv_{M,t})\sin(\frac{2\pi v_{M,t}}{\sqrt{k}})+(2Mtu_{M,t}+M^2v_{M,t}-t^2v_{M,t})\sinh(\frac{2\pi u_{M,t}}{\sqrt{k}}))}{2\sqrt{k}(M^2+t^2)^2\sqrt{M^2+t^2}(\sinh^2(\frac{\pi u_{M,t}}{\sqrt{k}}))+\sin^2(\frac{\pi v_{M,t}}{\sqrt{k}}))} .& \label{proofeq12} \tag{3.6} 
   \end{flalign*}

Let us consider the series $\sum\limits_{r=1}^{\infty} (-1)^re^{\pi i r \beta}H(r)$, and perform a Weierstrass M-test \cite[p.~103]{ComplexFunctionsReinhold}. By \eqref{proofeq12}, the series can be separately considered as a combination of two series, one corresponding to the part $\frac{it\pi^4}{30k^2(M^2+t^2)}- \frac{2\pi^2itM}{3k(M^2+t^2)^2} - it\frac{(3M^2-t^2)}{(M^2+t^2)^3}$ and the other corresponding to the remaining part. It is easily observed that the series corresponding to the part $\frac{it\pi^4}{30k^2(M^2+t^2)}- \frac{2\pi^2itM}{3k(M^2+t^2)^2} - it\frac{(3M^2-t^2)}{(M^2+t^2)^3}$ is absolutely and uniformly convergent for each $|\beta| \leq1$. Therefore we are left to deal with the the remaining series only. For each term in this remaining series we have

\begin{flalign*}
    &\left|(-1)^re^{\pi i r \beta}\frac{i\pi((r^2u_{r,t}-t^2u_{r,t}-2rtv_{r,t})\sin(\frac{2\pi v_{r,t}}{\sqrt{k}})+(2rtu_{r,t}+r^2v_{r,t}-t^2v_{r,t})\sinh(\frac{2\pi u_{r,t}}{\sqrt{k}}))}{2\sqrt{k}(r^2+t^2)^2\sqrt{r^2+t^2}(\sinh^2(\frac{\pi u_{r,t}}{\sqrt{k}})+\sin^2(\frac{\pi v_{r,t}}{\sqrt{k}}))}\right| &\\& \leq \frac{\pi}{2\sqrt{k}}\left|\frac{r^2u_{r,t}\sin(\frac{2\pi v_{r,t}}{\sqrt{k}})}{(r^2+t^2)^2\sqrt{r^2+t^2}(\sinh^2(\frac{\pi u_{r,t}}{\sqrt{k}})+\sin^2(\frac{\pi v_{r,t}}{\sqrt{k}}))} \right|+\frac{\pi t^2}{2\sqrt{k}}\left|\frac{u_{r,t} \sin(\frac{2\pi v_{r,t}}{\sqrt{k}})}{(r^2+t^2)^2\sqrt{r^2+t^2}(\sinh^2(\frac{\pi u_{r,t}}{\sqrt{k}})+\sin^2(\frac{\pi v_{r,t}}{\sqrt{k}}))} \right| &\\& +\frac{\pi |t|}{\sqrt{k}}\left|\frac{rv_{r,t} \sin(\frac{2\pi v_{r,t}}{\sqrt{k}})}{(r^2+t^2)^2\sqrt{r^2+t^2}(\sinh^2(\frac{\pi u_{r,t}}{\sqrt{k}})+\sin^2(\frac{\pi v_{r,t}}{\sqrt{k}}))} \right|+\frac{\pi |t|}{\sqrt{k}}\left|\frac{ru_{r,t} \sinh(\frac{2\pi u_{r,t}}{\sqrt{k}})}{(r^2+t^2)^2\sqrt{r^2+t^2}(\sinh^2(\frac{\pi u_{r,t}}{\sqrt{k}})+\sin^2(\frac{\pi v_{r,t}}{\sqrt{k}}))} \right|&\end{flalign*} \begin{flalign*}& + \frac{\pi}{2\sqrt{k}}\left|\frac{r^2v_{r,t}\sinh(\frac{2\pi u_{r,t}}{\sqrt{k}})}{(r^2+t^2)^2\sqrt{r^2+t^2}(\sinh^2(\frac{\pi u_{r,t}}{\sqrt{k}})+\sin^2(\frac{\pi v_{r,t}}{\sqrt{k}}))} \right|+\frac{\pi t^2}{2\sqrt{k}}\left|\frac{v_{r,t} \sinh(\frac{2\pi u_{r,t}}{\sqrt{k}})}{(r^2+t^2)^2\sqrt{r^2+t^2}(\sinh^2(\frac{\pi u_{r,t}}{\sqrt{k}})+\sin^2(\frac{\pi v_{r,t}}{\sqrt{k}}))} \right| &\\& < \frac{\pi}{2\sqrt{k}}\frac{r^2u_{r,t}|\sin(\frac{2\pi v_{r,t}}{\sqrt{k}})|}{r^5\sinh^2(\frac{\pi u_{r,t}}{\sqrt{k}}))} +\frac{\pi t^2}{2\sqrt{k}}\frac{u_{r,t} |\sin(\frac{2\pi v_{r,t}}{\sqrt{k}})|}{r^5\sinh^2(\frac{\pi u_{r,t}}{\sqrt{k}})} +\frac{\pi |t|}{\sqrt{k}}\frac{rv_{r,t} |\sin(\frac{2\pi v_{r,t}}{\sqrt{k}})|}{r^5\sinh^2(\frac{\pi u_{r,t}}{\sqrt{k}})} +\frac{\pi |t|}{\sqrt{k}}\frac{ru_{r,t} \sinh(\frac{2\pi u_{r,t}}{\sqrt{k}})}{r^5\sinh^2(\frac{\pi u_{r,t}}{\sqrt{k}})} &\\& + \frac{\pi}{2\sqrt{k}}\frac{r^2v_{r,t}\sinh(\frac{2\pi u_{r,t}}{\sqrt{k}})}{r^5\sinh^2(\frac{\pi u_{r,t}}{\sqrt{k}})} +\frac{\pi t^2}{2\sqrt{k}}\frac{v_{r,t} \sinh(\frac{2\pi u_{r,t}}{\sqrt{k}})}{r^5\sinh^2(\frac{\pi u_{r,t}}{\sqrt{k}})} 
    &\\& = \frac{\pi}{2\sqrt{k}}\frac{u_{r,t}|\sin(\frac{2\pi v_{r,t}}{\sqrt{k}})|}{r^3\sinh^2(\frac{\pi u_{r,t}}{\sqrt{k}}))} +\frac{\pi t^2}{2\sqrt{k}}\frac{u_{r,t} |\sin(\frac{2\pi v_{r,t}}{\sqrt{k}})|}{r^5\sinh^2(\frac{\pi u_{r,t}}{\sqrt{k}})} +\frac{\pi |t|}{\sqrt{k}}\frac{v_{r,t} |\sin(\frac{2\pi v_{r,t}}{\sqrt{k}})|}{r^4\sinh^2(\frac{\pi u_{r,t}}{\sqrt{k}})} +\frac{\pi |t|}{\sqrt{k}}\frac{u_{r,t} \sinh(\frac{2\pi u_{r,t}}{\sqrt{k}})}{r^4\sinh^2(\frac{\pi u_{r,t}}{\sqrt{k}})} &\\& + \frac{\pi}{2\sqrt{k}}\frac{v_{r,t}\sinh(\frac{2\pi u_{r,t}}{\sqrt{k}})}{r^3\sinh^2(\frac{\pi u_{r,t}}{\sqrt{k}})} +\frac{\pi t^2}{2\sqrt{k}}\frac{v_{r,t} \sinh(\frac{2\pi u_{r,t}}{\sqrt{k}})}{r^5\sinh^2(\frac{\pi u_{r,t}}{\sqrt{k}})}&\\& \leq \frac{\pi \eta_t}{2\sqrt{k}}\frac{|\sin(\frac{2\pi v_{r,t}}{\sqrt{k}})|}{r^{\frac{5}{2}}\sinh^2(\frac{\pi u_{r,t}}{\sqrt{k}}))} +\frac{\pi t^2 \eta_t}{2\sqrt{k}}\frac{ |\sin(\frac{2\pi v_{r,t}}{\sqrt{k}})|}{r^\frac{9}{2}\sinh^2(\frac{\pi u_{r,t}}{\sqrt{k}})} +\frac{\pi |t| \xi_t}{\sqrt{k}}\frac{ |\sin(\frac{2\pi v_{r,t}}{\sqrt{k}})|}{r^\frac{9}{2}\sinh^2(\frac{\pi u_{r,t}}{\sqrt{k}})} +\frac{\pi |t| \eta_t}{\sqrt{k}}\frac{ \sinh(\frac{2\pi u_{r,t}}{\sqrt{k}})}{r^\frac{7}{2}\sinh^2(\frac{\pi u_{r,t}}{\sqrt{k}})} &\\& + \frac{\pi\xi_t}{2\sqrt{k}}\frac{\sinh(\frac{2\pi u_{r,t}}{\sqrt{k}})}{r^\frac{7}{2}\sinh^2(\frac{\pi u_{r,t}}{\sqrt{k}})} +\frac{\pi t^2 \xi_t}{2\sqrt{k}}\frac{ \sinh(\frac{2\pi u_{r,t}}{\sqrt{k}})}{r^\frac{11}{2}\sinh^2(\frac{\pi u_{r,t}}{\sqrt{k}})} = M_r. & \label{proofeq13} \tag{3.7} 
\end{flalign*}

Where the last inequality follows by noting that $\frac{u_{r,t}}{\sqrt{r}}$ and $\sqrt{r}v_{r,t}$ are bounded above as functions of $r$ for $r\geq 1$ and letting $\frac{u_{r,t}}{\sqrt{r}} \leq \eta_t$ and $\sqrt{r}v_{r,t} \leq \xi_t$, where $\eta_t, \xi_t$ are non-negative and depend only on $t$. Now it is easily seen that $\sum\limits_{r=1}^{\infty} M_r$ converges. Thus by Weierstrass criterion, it follows that the series under investigation converges absolutely and uniformly for $|\beta|\leq 1$.

Next we consider the series $\sum\limits_{r=1}^{\infty} (-1)^re^{-\pi i r \beta}H(-r)$. We argue in the same manner as before so that we are only left to test the series 

\[
\sum\limits_{r=1}^{\infty} (-1)^re^{-\pi i r \beta}\frac{i\pi((r^2v_{r,t}-t^2v_{r,t}+2rtu_{r,t})\sin(\frac{2\pi u_{r,t}}{\sqrt{k}})+(-2rtv_{r,t}+r^2u_{r,t}-t^2u_{r,t})\sinh(\frac{2\pi v_{r,t}}{\sqrt{k}}))}{2\sqrt{k}(r^2+t^2)^2\sqrt{r^2+t^2}(\sinh^2(\frac{\pi v_{r,t}}{\sqrt{k}})+\sin^2(\frac{\pi u_{r,t}}{\sqrt{k}}))}.\label{proofeq14} \tag{3.8}
\]

For each term in this series we have

\begin{flalign*}
    & \left| (-1)^re^{-\pi i r \beta}\frac{i\pi((r^2v_{r,t}-t^2v_{r,t}+2rtu_{r,t})\sin(\frac{2\pi u_{r,t}}{\sqrt{k}})+(-2rtv_{r,t}+r^2u_{r,t}-t^2u_{r,t})\sinh(\frac{2\pi v_{r,t}}{\sqrt{k}}))}{2\sqrt{k}(r^2+t^2)^2\sqrt{r^2+t^2}(\sinh^2(\frac{\pi v_{r,t}}{\sqrt{k}})+\sin^2(\frac{\pi u_{r,t}}{\sqrt{k}}))} \right| &\\& \leq \frac{\pi}{2\sqrt{k}}\left|\frac{r^2v_{r,t}\sin(\frac{2\pi u_{r,t}}{\sqrt{k}})}{(r^2+t^2)^2\sqrt{r^2+t^2}(\sinh^2(\frac{\pi v_{r,t}}{\sqrt{k}})+\sin^2(\frac{\pi u_{r,t}}{\sqrt{k}}))} \right|+\frac{\pi t^2}{2\sqrt{k}}\left|\frac{v_{r,t} \sin(\frac{2\pi u_{r,t}}{\sqrt{k}})}{(r^2+t^2)^2\sqrt{r^2+t^2}(\sinh^2(\frac{\pi v_{r,t}}{\sqrt{k}}))+\sin^2(\frac{\pi u_{r,t}}{\sqrt{k}}))} \right| &\\& +\frac{\pi |t|}{\sqrt{k}}\left|\frac{ru_{r,t} \sin(\frac{2\pi u_{r,t}}{\sqrt{k}})}{(r^2+t^2)^2\sqrt{r^2+t^2}(\sinh^2(\frac{\pi v_{r,t}}{\sqrt{k}})+\sin^2(\frac{\pi u_{r,t}}{\sqrt{k}}))} \right|+\frac{\pi |t|}{\sqrt{k}}\left|\frac{rv_{r,t} \sinh(\frac{2\pi v_{r,t}}{\sqrt{k}})}{(r^2+t^2)^2\sqrt{r^2+t^2}(\sinh^2(\frac{\pi v_{r,t}}{\sqrt{k}})+\sin^2(\frac{\pi u_{r,t}}{\sqrt{k}}))} \right| &\\& + \frac{\pi}{2\sqrt{k}}\left|\frac{r^2u_{r,t}\sinh(\frac{2\pi v_{r,t}}{\sqrt{k}})}{(r^2+t^2)^2\sqrt{r^2+t^2}(\sinh^2(\frac{\pi v_{r,t}}{\sqrt{k}})+\sin^2(\frac{\pi u_{r,t}}{\sqrt{k}}))} \right| + \frac{\pi t^2}{2\sqrt{k}}\left|\frac{u_{r,t} \sinh(\frac{2\pi v_{r,t}}{\sqrt{k}})}{(r^2+t^2)^2\sqrt{r^2+t^2}(\sinh^2(\frac{\pi v_{r,t}}{\sqrt{k}})+\sin^2(\frac{\pi u_{r,t}}{\sqrt{k}}))} \right| &\\& < \frac{\pi}{2\sqrt{k}}\frac{r^2v_{r,t}|\sin(\frac{2\pi u_{r,t}}{\sqrt{k}})|}{r^5\sinh^2(\frac{\pi v_{r,t}}{\sqrt{k}}))} +\frac{\pi t^2}{2\sqrt{k}}\frac{v_{r,t} |\sin(\frac{2\pi u_{r,t}}{\sqrt{k}})|}{r^5\sinh^2(\frac{\pi v_{r,t}}{\sqrt{k}})} +\frac{\pi |t|}{\sqrt{k}}\frac{ru_{r,t} |\sin(\frac{2\pi u_{r,t}}{\sqrt{k}})|}{r^5\sinh^2(\frac{\pi v_{r,t}}{\sqrt{k}})} +\frac{\pi |t|}{\sqrt{k}}\frac{rv_{r,t} \sinh(\frac{2\pi v_{r,t}}{\sqrt{k}})}{r^5\sinh^2(\frac{\pi v_{r,t}}{\sqrt{k}})} &\\& + \frac{\pi}{2\sqrt{k}}\frac{r^2u_{r,t}\sinh(\frac{2\pi v_{r,t}}{\sqrt{k}})}{r^5\sinh^2(\frac{\pi v_{r,t}}{\sqrt{k}})} +\frac{\pi t^2}{2\sqrt{k}}\frac{u_{r,t} \sinh(\frac{2\pi v_{r,t}}{\sqrt{k}})}{r^5\sinh^2(\frac{\pi v_{r,t}}{\sqrt{k}})} &\end{flalign*} \begin{flalign*}& = \frac{\pi}{2\sqrt{k}}\frac{v_{r,t}|\sin(\frac{2\pi u_{r,t}}{\sqrt{k}})|}{r^3\sinh^2(\frac{\pi v_{r,t}}{\sqrt{k}}))} +\frac{\pi t^2}{2\sqrt{k}}\frac{v_{r,t} |\sin(\frac{2\pi u_{r,t}}{\sqrt{k}})|}{r^5\sinh^2(\frac{\pi v_{r,t}}{\sqrt{k}})} +\frac{\pi |t|}{\sqrt{k}}\frac{u_{r,t} |\sin(\frac{2\pi u_{r,t}}{\sqrt{k}})|}{r^4\sinh^2(\frac{\pi v_{r,t}}{\sqrt{k}})} +\frac{\pi |t|}{\sqrt{k}}\frac{v_{r,t} \sinh(\frac{2\pi v_{r,t}}{\sqrt{k}})}{r^4\sinh^2(\frac{\pi v_{r,t}}{\sqrt{k}})} &\\& + \frac{\pi}{2\sqrt{k}}\frac{u_{r,t}\sinh(\frac{2\pi v_{r,t}}{\sqrt{k}})}{r^3\sinh^2(\frac{\pi v_{r,t}}{\sqrt{k}})} +\frac{\pi t^2}{2\sqrt{k}}\frac{u_{r,t} \sinh(\frac{2\pi v_{r,t}}{\sqrt{k}})}{r^5\sinh^2(\frac{\pi v_{r,t}}{\sqrt{k}})} &\\& < \frac{\pi}{2\sqrt{k}}\frac{v_{r,t}}{r^3\sinh^2(\frac{\pi v_{r,t}}{\sqrt{k}}))} +\frac{\pi t^2}{2\sqrt{k}}\frac{v_{r,t} }{r^5\sinh^2(\frac{\pi v_{r,t}}{\sqrt{k}})} +\frac{\pi |t|}{\sqrt{k}}\frac{u_{r,t} }{r^4\sinh^2(\frac{\pi v_{r,t}}{\sqrt{k}})} +\frac{\pi |t|}{\sqrt{k}}\frac{v_{r,t} \sinh(\frac{2\pi v_{r,t}}{\sqrt{k}})}{r^4\sinh^2(\frac{\pi v_{r,t}}{\sqrt{k}})} &\\& + \frac{\pi}{2\sqrt{k}}\frac{u_{r,t}\sinh(\frac{2\pi v_{r,t}}{\sqrt{k}})}{r^3\sinh^2(\frac{\pi v_{r,t}}{\sqrt{k}})} +\frac{\pi t^2}{2\sqrt{k}}\frac{u_{r,t} \sinh(\frac{2\pi v_{r,t}}{\sqrt{k}})}{r^5\sinh^2(\frac{\pi v_{r,t}}{\sqrt{k}})}.
    \label{proofeq15} \tag{3.9}
\end{flalign*}

Before proceeding further, we note that as a function of $r$, $\sqrt{r}v_{r,t}$ is bounded below for $r\geq 1$, thus we let $\sqrt{r}v_{r,t} \geq \sigma_t $, where $\sigma_t$ is non-negative and depends only on $t$. We also note that as functions of $r$, $\frac{\pi v_{r,t}}{\sqrt{k}\sinh(\frac{\pi v_{r,t}}{\sqrt{k}})}$ and $\frac{\sinh(\frac{2\pi v_{r,t}}{\sqrt{k}})}{\sinh(\frac{\pi v_{r,t}}{\sqrt{k}})}$ are bounded above for $r \geq 1$. Thus we let $\frac{\pi v_{r,t}}{\sqrt{k}\sinh(\frac{\pi v_{r,t}}{\sqrt{k}})} \leq \Xi$ and $\frac{\sinh(\frac{2\pi v_{r,t}}{\sqrt{k}})}{\sinh(\frac{\pi v_{r,t}}{\sqrt{k}})} \leq \Theta$, for $r \geq 1$ where $\Xi$ and $\Theta$ are non-negative and independent of $r$. Thus, we have that 

\begin{flalign*}
    & \left| (-1)^re^{-\pi i r \beta}\frac{i\pi((r^2v_{r,t}-t^2v_{r,t}+2rtu_{r,t})\sin(\frac{2\pi u_{r,t}}{\sqrt{k}})+(-2rtv_{r,t}+r^2u_{r,t}-t^2u_{r,t})\sinh(\frac{2\pi v_{r,t}}{\sqrt{k}}))}{2\sqrt{k}(r^2+t^2)^2\sqrt{r^2+t^2}(\sinh^2(\frac{\pi v_{r,t}}{\sqrt{k}})+\sin^2(\frac{\pi u_{r,t}}{\sqrt{k}}))} \right| &\\& < \frac{\Xi^2\sqrt{k}}{2\sigma_t\pi}\frac{1}{r^\frac{5}{2}} +\frac{\Xi^2\sqrt{k} t^2}{2\sigma_t\pi}\frac{1 }{r^\frac{9}{2}} +\frac{\Xi^2\eta_t\sqrt{k} |t|}{\sigma_t^2\pi}\frac{1 }{r^{\frac{5}{2}}} +\Xi \Theta |t|\frac{1}{r^4} + \frac{\Xi \Theta \eta_t}{2\sigma_t}\frac{1}{r^{2}} +\frac{\Xi \Theta \eta_t t^2}{2\sigma_t}\frac{1}{r^{4}} = M'_r. \label{proofeq16} \tag{3.10}
\end{flalign*}

Now it is obvious that the series $\sum\limits_{r=1}^{\infty} M'_r $ converges and thus, by Weierstrass criterion we conclude that the series under investigation converges absolutely and uniformly for $|\beta|\leq 1$, proving the lemma.  
\end{proof}

Therefore, by \textbf{Lemma 6}, \textbf{Lemma 7} and \textbf{Theorem 5}, we are now able to express $q_k(N)$ in an analytic expression in the following manner. 

\begin{flalign*}
        &   \frac{q_k(N)}{N^2}  &\\& = \frac{\sinh(\pi t)}{2i\pi }\sum\limits_{r=1}^{\infty} (-1)^{r}\left(\frac{it\pi^4}{30k^2((r-N)^2+t^2)}- \frac{2\pi^2it(r-N)}{3k((r-N)^2+t^2)^2} - it\frac{(3(r-N)^2-t^2)}{((r-N)^2+t^2)^3} \right. &\\& \left.+\frac{it\pi^4}{30k^2((r+N)^2+t^2)}+ \frac{2\pi^2it(r+N)}{3k((r+N)^2+t^2)^2} - it\frac{(3(r+N)^2-t^2)}{((r+N)^2+t^2)^3} +\frac{i\pi}{2\sqrt{k}} (G_{r-N,t,k}+G_{-r-N,t,k})\right)\int_{0}^{1} \frac{\cos(\pi r \beta)}{\cosh(\pi\beta t)} d\beta &\\& +\left( \frac{it\pi^4}{30k^2(N^2+t^2)}+ \frac{2\pi^2itN}{3k(N^2+t^2)^2} - it\frac{(3N^2-t^2)}{(N^2+t^2)^3}+\frac{i\pi}{2\sqrt{k}} G_{-N,t,k}\right)\frac{\sinh(\pi t)\arctan(\tanh(\frac{\pi t}{2}))}{i\pi^2t}  . \label{proofeq17} \tag{3.11}
    \end{flalign*}

Where $G_{M,t,k}$ is just the expression defined in \eqref{Gdef}. Making further manipulations we have 

\begin{flalign*}
        &   \frac{q_k(N)}{N^2}  &\\& = \frac{(-1)^Nt \pi^3\sinh(\pi t)}{60k^2}\int_{0}^{1}\frac{\cos(\pi N \beta)}{\cosh(\pi \beta t)}\left( \sum\limits_{r=-\infty}^{\infty} \frac{(-1)^r\cos(\pi r \beta)}{r^2+t^2}\right) d\beta &\\& + \frac{(-1)^N\pi t\sinh(\pi t)}{3k}\int_{0}^{1}\frac{\sin(\pi N \beta)}{\cosh(\pi \beta t)}\left( \sum\limits_{r=-\infty}^{\infty} \frac{(-1)^r r\sin(\pi r \beta)}{(r^2+t^2)^2}\right) d\beta &\\& +(-1)^{N-1} \frac{t \sinh(\pi t)}{2\pi}\int_{0}^{1}\frac{\cos(\pi N \beta)}{\cosh(\pi \beta t)}\left( \sum\limits_{r=-\infty}^{\infty} \frac{(-1)^r(3r^2-t^2)\cos(\pi r \beta)}{(r^2+t^2)^3}\right) d\beta 
        + \frac{\sinh(\pi t) \arctan(\tanh(\frac{\pi t}{2}))}{2\pi t\sqrt{k}}G_{-N,t,k} &\\& +\frac{\sinh(\pi t)}{4\sqrt{k}}\sum\limits_{r=1}^{\infty}(-1)^r(G_{r-N,t,k}+G_{-r-N,t,k})\int_{0}^{1} \frac{\cos(\pi r \beta)}{\cosh(\pi\beta t)} d\beta .& \label{proofeq18} \tag{3.12}
    \end{flalign*}

Now since for $|\beta|\leq 1$, $\sum\limits_{r=-\infty}^{\infty} \frac{(-1)^r \cos(\pi r \beta)}{r^2+t^2}=\frac{\pi \cosh(\pi \beta t)}{t \sinh(\pi t)}$, $\sum\limits_{r=-\infty}^{\infty} \frac{(-1)^r r \sin(\pi r \beta)}{(r^2+t^2)^2}=\frac{\pi^2}{2t}\left( \frac{\beta \cosh(\pi \beta t)}{\sinh(\pi t)} - \frac{\sinh(\pi \beta t) \cosh(\pi t)}{\sinh^2(\pi t)}\right)$ and $\sum\limits_{r=-\infty}^{\infty} \frac{(-1)^r(3r^2-t^2)\cos(\pi r \beta)}{(r^2+t^2)^3} = -\frac{\pi^3 \beta^2 \cosh(\pi \beta t)}{2t \sinh(\pi t)} + \frac{\pi^2 \beta \sinh(\pi \beta t) \cosh(\pi t)}{t \sinh^2(\pi t)}- \frac{\pi^3 \cosh(\pi \beta t) \cosh^2(\pi t)}{t \sinh^3(\pi t)}+\frac{\pi^3 \cosh(\pi \beta t) \sinh(\pi t)}{2t \sinh^2(\pi t)}$, direct evaluation of the expression in the right hand side yields

\begin{flalign*}
        &  \frac{(1+(-1)^{N-1})\pi^2 \coth(\pi t) }{6Nk}-\frac{(1+(-1)^{N-1}) \coth(\pi t)}{2\pi N^2} + \frac{(1+(-1)^{N-1})}{2N^2} - \frac{\pi^2}{6Nk }&\\& + \frac{(-1)^N \pi^3 \coth(\pi t)}{3k} \int_{0}^{1} \frac{\sin(\pi N \beta)}{e^{2\pi \beta t}+1} d\beta + (-1)^N \pi \coth(\pi t) \int_{0}^{1} \frac{\beta \cos(\pi N \beta)}{e^{2\pi \beta t}+1} d\beta + \frac{\sinh(\pi t) \arctan(\tanh(\frac{\pi t}{2}))}{2\pi t\sqrt{k}}G_{-N,t,k} &\\& +\frac{\sinh(\pi t)}{4\sqrt{k}}\sum\limits_{r=1}^{\infty}(-1)^r(G_{r-N,t,k}+G_{-r-N,t,k})\int_{0}^{1} \frac{\cos(\pi r \beta)}{\cosh(\pi\beta t)} d\beta .
        \label{proofeq19} \tag{3.13}
    \end{flalign*}

Therefore we state the following proposition giving an analytic representation to $q_k(N)$.
\vspace{0.7cm}
\begin{proposition}
    Let $G_{M,t,k}$ be defined as in \eqref{Gdef}. Then, for each natural number $N$, and real $t \neq 0$

    \begin{flalign*}
        &   \frac{q_k(N)}{N^2}  &\\& = (1+(-1)^{N-1})\left\{\frac{\pi^2 \coth(\pi t) }{6Nk}-\frac{ \coth(\pi t)}{2\pi N^2} + \frac{1}{2N^2} \right\} - \frac{\pi^2}{6Nk }&\\& + \frac{(-1)^N \pi^3 \coth(\pi t)}{3k} \int_{0}^{1} \frac{\sin(\pi N \beta)}{e^{2\pi \beta t}+1} d\beta + (-1)^N \pi \coth(\pi t) \int_{0}^{1} \frac{\beta \cos(\pi N \beta)}{e^{2\pi \beta t}+1} d\beta + \frac{\sinh(\pi t) \arctan(\tanh(\frac{\pi t}{2}))}{2\pi t\sqrt{k}}G_{-N,t,k} &\\& +\frac{\sinh(\pi t)}{4\sqrt{k}}\sum\limits_{r=1}^{\infty}(-1)^r(G_{r-N,t,k}+G_{-r-N,t,k})\int_{0}^{1} \frac{\cos(\pi r \beta)}{\cosh(\pi\beta t)} d\beta .&
         \label{analyticFormulaForq} \tag{3.14}
    \end{flalign*}

    And for $N$ either zero or a negative integer, and real $t \neq 0$

    \begin{flalign*}
        &   0 = \mu_N + \frac{(-1)^N \pi^3 \coth(\pi t)}{3k} \int_{0}^{1} \frac{\sin(\pi N \beta)}{e^{2\pi \beta t}+1} d\beta + (-1)^N \pi \coth(\pi t) \int_{0}^{1} \frac{\beta \cos(\pi N \beta)}{e^{2\pi \beta t}+1} d\beta + \frac{\sinh(\pi t) \arctan(\tanh(\frac{\pi t}{2}))}{2\pi t\sqrt{k}}G_{-N,t,k} &\\& +\frac{\sinh(\pi t)}{4\sqrt{k}}\sum\limits_{r=1}^{\infty}(-1)^r(G_{r-N,t,k}+G_{-r-N,t,k})\int_{0}^{1} \frac{\cos(\pi r \beta)}{\cosh(\pi\beta t)} d\beta ,&
         \label{analyticFormulaFor0} \tag{3.15}
    \end{flalign*}

    where 
    \[
    \mu_N = \begin{cases}
        \frac{ \pi^2}{12} - \frac{\pi \coth(\pi t)}{4} & \textit{, if } N = 0, \\
        (1+(-1)^{N-1})\left\{\frac{\pi^2 \coth(\pi t) }{6Nk}-\frac{ \coth(\pi t)}{2\pi N^2} + \frac{1}{2N^2} \right\} - \frac{\pi^2}{6Nk } & \textit{, if else.}
    \end{cases}
    \]
\end{proposition}

Now we proceed to argue that such an analytic representation is always possible for all $q_{k,s}(N)$. We state the following lemma, which is obvious from the definition of $q_{k,s}(N)$.
\vspace{0.7cm}
\begin{lemma} For natural numbers $k$ and $s$ and for arbitrary complex value $z$, we have
    \[
    \sum\limits_{n=1}^{\infty} \frac{q_{k,s}(n)}{n^2(n+z)} = \sum\limits_{n=1}^{\infty} \frac{1}{k^2n^{4s}(kn^{2s}+z)}.
         \label{generatingSeriesForqks} \tag{3.16}
    \]
\end{lemma}

Then, we argue that for $F(z) = \sum\limits_{n=1}^{\infty} \frac{1}{k^2n^{4s}(kn^{2s}+z)}$, the series $\sum\limits_{r=1}^{\infty}(-1)^re^{\pi i r \beta}(F(r-it)-F(r+it))$ and $\sum\limits_{r=1}^{\infty} (-1)^re^{-\pi i r \beta}(F(-r-it)-F(-r+it))$ converges absolutely and uniformly for $|\beta| \leq 1$. Thus, we shall prove that
\vspace{0.7cm}
\begin{lemma} For natural numbers $k$ and $s$, if $F(z)$ is given by
\[F(z) = \sum\limits_{n=1}^{\infty} \frac{1}{k^2n^{4s}(kn^{2s}+z)},\label{2generatingSeriesForqks} \tag{3.17}\]

then the series $\sum\limits_{r=1}^{\infty} (-1)^re^{\pi i r \beta}(F(r-it)-F(r+it))$ and $\sum\limits_{r=1}^{\infty} (-1)^re^{-\pi i r \beta}(F(-r-it)-F(-r+it))$ converges absolutely and uniformly for all natural numbers $k$ and $s$, and $|\beta| \leq 1$.
    
\end{lemma}

\begin{proof}[Proof of Lemma 10]
    Let us compute $F(M-it)-F(M+it)$
\begin{flalign*}
    & F(M-it)-F(M+it) &\\& = \sum\limits_{n=1}^{\infty} \frac{1}{k^2n^{4s}(kn^{2s}+M-it)} - \sum\limits_{n=1}^{\infty} \frac{1}{k^2n^{4s}(kn^{2s}+M+it)} &\\& = \sum\limits_{n=1}^{\infty} \frac{2it}{k^2n^{4s}((kn^{2s}+M)^2+t^2)}& \label{proofeq20} \tag{3.18}
\end{flalign*}

Therefore, for each term in the series $\sum\limits_{r=1}^{\infty}(-1)^re^{\pi i r \beta}(F(r-it)-F(r+it))$ and $\sum\limits_{r=1}^{\infty}(-1)^re^{-\pi i r \beta}(F(-r-it)-F(-r+it))$, we have

\begin{flalign*}
& \left| (-1)^re^{\pi i r \beta} \sum\limits_{n=1}^{\infty} \frac{2it}{k^2n^{4s}((kn^{2s}+r)^2+t^2)}\right| &\\& \leq  2 |t|\sum\limits_{n=1}^{\infty} \frac{1}{k^2n^{4s}((kn^{2s}+r)^2+t^2)}  &\\& \leq 2 |t|\sum\limits_{n=1}^{\infty} \frac{1}{k^2n^{4}((kn^{2}+r)^2+t^2)} = |H(r)|, & \label{proofeq21} \tag{3.19}
\end{flalign*}

and

\begin{flalign*}
& \left| (-1)^re^{-\pi i r \beta} \sum\limits_{n=1}^{\infty} \frac{2it}{k^2n^{4s}((kn^{2s}-r)^2+t^2)}\right| &\\& \leq  2 |t|\sum\limits_{n=1}^{\infty} \frac{1}{k^2n^{4s}((kn^{2s}-r)^2+t^2)}  &\\& \leq 2 |t|\sum\limits_{n=1}^{\infty} \frac{1}{k^2n^{4}((kn^{2}-r)^2+t^2)} = |H(-r)|. & \label{proofeq22} \tag{3.20}
\end{flalign*}

Where $H$ is given by \eqref{proofeq12}. And since by \textbf{Lemma 7}, $\sum\limits_{r=1}^{\infty} |H(r)|$ and $\sum\limits_{r=1}^{\infty} |H(-r)|$ converges, by Weierstrass criterion we must have $\sum\limits_{r=1}^{\infty} (-1)^re^{\pi i r \beta}(F(r-it)-F(r+it))$ and $\sum\limits_{r=1}^{\infty} (-1)^re^{-\pi i r \beta}(F(-r-it)-F(-r+it))$ converges absolutely and uniformly for all natural numbers $k$ and $s$, and $|\beta| \leq 1$.

\end{proof}

Therefore, in view of \textbf{Lemma 9}, \textbf{Lemma 10} and \textbf{Theorem 5}, we now present an analytic expression for $q_{k,s}(N)$.
\begin{proposition} Let $k$ and $s$ be natural numbers, and let $H_{k,s}(z) = \sum\limits_{n=1}^{\infty} \frac{1}{k^2n^{4s}((kn^{2s}+z)^2+t^2)}$, then for each natural number $N$ and real $t \neq 0$, we have 
    \begin{flalign*}
        &   \frac{q_{k,s}(N)}{N^2}  =  \frac{t\sinh(\pi t)}{\pi }\sum\limits_{r=1}^{\infty} (-1)^{r}(H_{k,s}(r-N)+H_{k,s}(-r-N))\int_{0}^{1} \frac{\cos(\pi r \beta)}{\cosh(\pi\beta t)} d\beta +\frac{2H(-N)\sinh(\pi t)\arctan(\tanh(\frac{\pi t}{2}))}{\pi^2}  .& & \label{analyticFormulaForqks} \tag{3.21}
    \end{flalign*}
The expression in the right hand side vanishes when $N$ is either zero or a negative integer.
\end{proposition}

Lastly we remark that in the formulae \eqref{analyticFormulaForq} and \eqref{analyticFormulaFor0}, the part 

\[
\frac{\sinh(\pi t) \arctan(\tanh(\frac{\pi t}{2}))}{2\pi t\sqrt{k}}G_{-N,t,k}  +\frac{\sinh(\pi t)}{4\sqrt{k}}\sum\limits_{r=1}^{\infty}(-1)^r(G_{r-N,t,k}+G_{-r-N,t,k})\int_{0}^{1} \frac{\cos(\pi r \beta)}{\cosh(\pi\beta t)} d\beta,
\]

can be viewed as a bilateral series over $r$, as follows 

\[
\frac{\sinh(\pi t)}{4\sqrt{k}}\sum\limits_{r=-\infty}^{\infty}(-1)^rG_{r-N,t,k}\int_{0}^{1} \frac{\cos(\pi r \beta)}{\cosh(\pi\beta t)} d\beta.
\]

So that for an integer $c$ such that $N+c > 0$, we have

\begin{flalign*}
        &   \frac{q_k(N+c)}{(N+c)^2}  &\\& = (1+(-1)^{N+c-1})\left\{\frac{\pi^2 \coth(\pi t) }{6(N+c)k}-\frac{ \coth(\pi t)}{2\pi (N+c)^2} + \frac{1}{2(N+c)^2} \right\} - \frac{\pi^2}{6(N+c)k }&\\& + \frac{(-1)^{N+c} \pi^3 \coth(\pi t)}{3k} \int_{0}^{1} \frac{\sin(\pi (N+c) \beta)}{e^{2\pi \beta t}+1} d\beta + (-1)^{N+c} \pi \coth(\pi t) \int_{0}^{1} \frac{\beta \cos(\pi (N+c) \beta)}{e^{2\pi \beta t}+1} d\beta &\\& + \frac{\sinh(\pi t)}{4\sqrt{k}}\sum\limits_{r=-\infty}^{\infty}(-1)^rG_{r-N-c,t,k}\int_{0}^{1} \frac{\cos(\pi r \beta)}{\cosh(\pi\beta t)} d\beta  &\\& = (1+(-1)^{N+c-1})\left\{\frac{\pi^2 \coth(\pi t) }{6(N+c)k}-\frac{ \coth(\pi t)}{2\pi (N+c)^2} + \frac{1}{2(N+c)^2} \right\} - \frac{\pi^2}{6(N+c)k }&\\& + \frac{(-1)^{N+c} \pi^3 \coth(\pi t)}{3k} \int_{0}^{1} \frac{\sin(\pi (N+c) \beta)}{e^{2\pi \beta t}+1} d\beta + (-1)^{N+c} \pi \coth(\pi t) \int_{0}^{1} \frac{\beta \cos(\pi (N+c) \beta)}{e^{2\pi \beta t}+1} d\beta &\\& + \frac{\sinh(\pi t)G_{-N,t,k}}{4\sqrt{k}}\int_{0}^{1} \frac{\cos(\pi c\beta)}{\cosh(\pi\beta t)} d\beta +\frac{\sinh(\pi t)}{4\sqrt{k}}\sum\limits_{r=1}^{\infty}(-1)^rG_{r-N,t,k}\int_{0}^{1} \frac{\cos(\pi (r +c)\beta)}{\cosh(\pi\beta t)} d\beta &\\& +\frac{\sinh(\pi t)}{4\sqrt{k}}\sum\limits_{r=1}^{\infty}(-1)^rG_{-r-N,t,k}\int_{0}^{1} \frac{\cos(\pi (r -c)\beta)}{\cosh(\pi\beta t)} d\beta. \label{proofeq23} \tag{3.22}
    \end{flalign*}

    And the counterpart for $N+c \leq 0$ is 
    \begin{flalign*}
        &   0 = \mu_{N+c} + \frac{(-1)^{N+c} \pi^3 \coth(\pi t)}{3k} \int_{0}^{1} \frac{\sin(\pi (N+c) \beta)}{e^{2\pi \beta t}+1} d\beta + (-1)^{N+c} \pi \coth(\pi t) \int_{0}^{1} \frac{\beta \cos(\pi (N+c) \beta)}{e^{2\pi \beta t}+1} d\beta &\\& + \frac{\sinh(\pi t)G_{-N,t,k}}{4\sqrt{k}}\int_{0}^{1} \frac{\cos(\pi c\beta)}{\cosh(\pi\beta t)} d\beta +\frac{\sinh(\pi t)}{4\sqrt{k}}\sum\limits_{r=1}^{\infty}(-1)^rG_{r-N,t,k}\int_{0}^{1} \frac{\cos(\pi (r +c)\beta)}{\cosh(\pi\beta t)} d\beta &\\& +\frac{\sinh(\pi t)}{4\sqrt{k}}\sum\limits_{r=1}^{\infty}(-1)^rG_{-r-N,t,k}\int_{0}^{1} \frac{\cos(\pi (r -c)\beta)}{\cosh(\pi\beta t)} d\beta , \label{proofeq24} \tag{3.23}
    \end{flalign*}

    where 
    \[
    \mu_R = \begin{cases}
        \frac{ \pi^2}{12} - \frac{\pi \coth(\pi t)}{4} & \textit{, if } R = 0, \\
        (1+(-1)^{R-1})\left\{\frac{\pi^2 \coth(\pi t) }{6Rk}-\frac{ \coth(\pi t)}{2\pi R^2} + \frac{1}{2R^2} \right\} - \frac{\pi^2}{6Rk } & \textit{, if else.}
    \end{cases}
    \]

We now proceed to evaluate in the next section, the following three types of integrals appearing in our formula for $q_k(N)$ 

\begin{flalign*}
    I_{q,t} &= \int_{0}^{1} \frac{\sin(\pi q \beta)}{e^{2\pi \beta t}+1} d\beta, \label{proofeq25} \tag{3.24}\\ 
    K_{q,t} &= \int_{0}^{1} \frac{\beta \cos(\pi q \beta)}{e^{2\pi \beta t}+1} d\beta, \label{proofeq26} \tag{3.25}\\
    J_{q,t} &= \int_{0}^{1} \frac{\cos(\pi q \beta)}{\cosh(\pi\beta t)} d\beta.\label{proofeq27} \tag{3.26}
\end{flalign*}

\section{Evaluating $I_{q,t}$, $K_{q,t}$ and $J_{q,t}$}

Here we shall seek to make some improvements on the pair of equalities \eqref{proofeq23} and \eqref{proofeq24}. The following lemmas are certainly not new, yet it would be more preferable to derive them on our own rather than delving into vast mathematical literature to search for an equivalent. 
\vspace{0.7cm}
\begin{lemma} For an integer $q$ and real $t > 0$, we have

\begin{flalign*}
    I_{q,t} &= i_q + \frac{(-1)^{q-1}q}{\pi} \sum\limits_{r=1}^{\infty} \frac{(-1)^{r-1}e^{-2\pi t r}}{4t^2r^2+q^2}, \label{proofeq28} \tag{4.1}
     \\
    K_{q,t} &= k_q + \frac{(-1)^{q-1}t}{\pi} \sum\limits_{r=1}^{\infty} \frac{(-1)^{r-1}re^{-2\pi t r}}{4t^2r^2+q^2} + \frac{(-1)^{q-1}}{\pi^2} \sum\limits_{r=1}^{\infty} \frac{(-1)^{r-1}e^{-2\pi t r}(4t^2r^2-q^2)}{(4t^2r^2+q^2)^2},\label{proofeq29} \tag{4.2}
    \\
    J_{q,t} &= \frac{1}{2t \cosh(\frac{\pi q}{2t})}+ \frac{(-1)^{q-1}2t}{\pi}\sum\limits_{r=0}^{\infty}\frac{(-1)^{r}(2r+1)e^{-\pi t(2r+1)}}{t^2(2r+1)^2+q^2}. \label{proofeq30} \tag{4.3}
\end{flalign*}

Where 

\begin{flalign*}
    i_q = &\begin{cases} 0 &, \textit{ if } q = 0, \\
    \frac{1}{2\pi q}- \frac{1}{4t\sinh(\frac{\pi q}{2t})} &, \textit{ if else,} 
    \end{cases}
    \\
    k_q = &\begin{cases} \frac{1}{48t^2} &, \textit{ if } q = 0, \\
     \frac{1}{2\pi^2q^2} - \frac{t^2}{q^2 \pi^2}- \frac{1}{4\pi q t \sinh(\frac{\pi q}{2t})} + \frac{t}{4\pi q \sinh(\frac{\pi q}{2t})} + \frac{\cosh(\frac{\pi q}{2t})}{8\sinh^2(\frac{\pi q}{2t})} &, \textit{ if else.} 
    \end{cases}
\end{flalign*}
    
\end{lemma}

\begin{proof}[Proof of Lemma 12]
    In each of the integrals, we expand the parts $\frac{1}{e^{2\pi \beta t}+1}= \frac{e^{-2\pi \beta t}}{e^{-2\pi \beta t}+1}$ and $\frac{1}{\cosh(\pi \beta t)} = \frac{2e^{-\pi \beta t}}{e^{-2\pi \beta t}+1}$ as geometric series and integrate term-wise. Elementary manipulations and application of the formulae 

    \begin{flalign*}
        \sum\limits_{r=1}^{\infty} \frac{(-1)^{r-1}}{r^2+z^2} &= \frac{1}{2z^2}-\frac{\pi}{2z\sinh(\pi z)},\label{proofeq31} \tag{4.4} \\
        \sum\limits_{r=1}^{\infty} \frac{(-1)^{r-1}}{(r^2+z^2)^2} &= \frac{1}{2z^4}- \frac{\pi}{4z^3\sinh(\pi x)} - \frac{\pi^2 \cosh(\pi z)}{4z^2 \sinh^2(\pi z)},\label{proofeq32} \tag{4.5} \\ 
        \sum\limits_{r=0}^{\infty} \frac{(-1)^r(2r+1)}{(2r+1)^2+z} &= \frac{\pi}{4 \cosh(\frac{\pi z}{2})}, \label{proofeq33} \tag{4.6}
    \end{flalign*}
    yields the desired equalities.
\end{proof}

We can now improve the statement of \eqref{proofeq23} and \eqref{proofeq24} as follows. For $N+c > 0$

\begin{flalign*}
        &   \frac{q_k(N+c)}{(N+c)^2}  &\\& = (1+(-1)^{N+c-1})\left\{\frac{\pi^2 \coth(\pi t) }{6(N+c)k}-\frac{ \coth(\pi t)}{2\pi (N+c)^2} + \frac{1}{2(N+c)^2} \right\} - \frac{\pi^2}{6(N+c)k }&\\& + \frac{(-1)^{N+c} \pi^3 \coth(\pi t)}{3k} \left(\frac{1}{2\pi (N+c)}- \frac{1}{4t\sinh(\frac{\pi (N+c)}{2t})} + \frac{(-1)^{N+c-1}(N+c)}{\pi} \sum\limits_{r=1}^{\infty} \frac{(-1)^{r-1}e^{-2\pi t r}}{4t^2r^2+(N+c)^2}\right) &\\&+ (-1)^{N+c} \pi \coth(\pi t) \left( \frac{1}{2\pi^2(N+c)^2} - \frac{t^2}{(N+c)^2 \pi^2}- \frac{1}{4\pi (N+c) t \sinh(\frac{\pi (N+c)}{2t})} + \frac{t}{4\pi (N+c) \sinh(\frac{\pi (N+c)}{2t})} \right. &\\& \left. + \frac{\cosh(\frac{\pi (N+c)}{2t})}{8\sinh^2(\frac{\pi (N+c)}{2t})} +\frac{(-1)^{N+c-1}t}{\pi} \sum\limits_{r=1}^{\infty} \frac{(-1)^{r-1}re^{-2\pi t r}}{4t^2r^2+(N+c)^2} + \frac{(-1)^{N+c-1}}{\pi^2} \sum\limits_{r=1}^{\infty} \frac{(-1)^{r-1}e^{-2\pi t r}(4t^2r^2-(N+c)^2)}{(4t^2r^2+(N+c)^2)^2}\right) &\\& + \frac{\sinh(\pi t)G_{-N,t,k}}{4\sqrt{k}}\left(\frac{1}{2t \cosh(\frac{\pi c}{2t})}+ \frac{(-1)^{c-1}2t}{\pi}\sum\limits_{r=0}^{\infty}\frac{(-1)^{r}(2r+1)e^{-\pi t(2r+1)}}{t^2(2r+1)^2+c^2} \right) 
        &\end{flalign*} \begin{flalign*}&+\frac{\sinh(\pi t)}{4\sqrt{k}}\sum\limits_{r=1}^{\infty}(-1)^rG_{r-N,t,k}\left(\frac{1}{2t \cosh(\frac{\pi (r+c)}{2t})}+ \frac{(-1)^{r+c-1}2t}{\pi}\sum\limits_{m=0}^{\infty}\frac{(-1)^{m}(2m+1)e^{-\pi t(2m+1)}}{t^2(2m+1)^2+(r+c)^2} \right) &\\& +\frac{\sinh(\pi t)}{4\sqrt{k}}\sum\limits_{r=1}^{\infty}(-1)^rG_{-r-N,t,k}\left(\frac{1}{2t \cosh(\frac{\pi (r-c)}{2t})}+ \frac{(-1)^{r-c-1}2t}{\pi}\sum\limits_{m=0}^{\infty}\frac{(-1)^{m}(2m+1)e^{-\pi t(2m+1)}}{t^2(2m+1)^2+(r-c)^2} \right)
        &\\& = \frac{\pi^2}{3(N+c)k(e^{2\pi t}-1)}+ \frac{\pi - \coth(\pi t)}{2\pi (N+c)^2}+ \frac{(-1)^{N+c}(1-t^2)\coth(\pi t)}{\pi (N+c)^2} + \frac{(-1)^{N+c-1} \pi^3 \coth(\pi t)}{12kt \sinh(\frac{\pi (N+c)}{2t})} &\\& + \frac{(-1)^{N+c}\coth(\pi t)(t-\frac{1}{t})}{4(N+c)\sinh(\frac{\pi (N+c)}{2t})}+ \frac{(-1)^{N+c}\pi \cosh(\frac{\pi(N+c)}{2t})\coth(\pi t)}{8\sinh^2(\frac{\pi (N+c)}{2t})}- \frac{(N+c)\pi^2\coth(\pi t)}{3k} \sum\limits_{r=1}^{\infty} \frac{(-1)^{r-1}e^{-2\pi t r}}{4t^2r^2+(N+c)^2} &\\& - t\coth(\pi t) \sum\limits_{r=1}^{\infty} \frac{(-1)^{r-1}re^{-2\pi t r}}{4t^2r^2+(N+c)^2} - \frac{\coth(\pi t)}{\pi}\sum\limits_{r=1}^{\infty} \frac{(-1)^{r-1}e^{-2\pi t r}(4t^2r^2-(N+c)^2)}{(4t^2r^2+(N+c)^2)^2} + \frac{\sinh(\pi t) G_{-N,t,k}}{8\sqrt{k} t \cosh(\frac{\pi c}{2t})} &\\& + \frac{(-1)^{c-1}t\sinh(\pi t)G_{-N,t,k}}{2\sqrt{k}\pi}\sum\limits_{r=0}^{\infty}\frac{(-1)^{r}(2r+1)e^{-\pi t(2r+1)}}{t^2(2r+1)^2+c^2}
        + \frac{\sinh(\pi t)}{8\sqrt{k}t}\sum\limits_{r=1}^{\infty}(-1)^r\frac{G_{r-N,t,k}}{ \cosh(\frac{\pi (r+c)}{2t})} &\\&+\frac{\sinh(\pi t)}{8\sqrt{k}t}\sum\limits_{r=1}^{\infty}(-1)^r\frac{G_{-r-N,t,k}}{ \cosh(\frac{\pi (r-c)}{2t})}  +\frac{(-1)^{c-1}t\sinh(\pi t)}{2\sqrt{k} \pi}\sum\limits_{r=1}^{\infty} \sum\limits_{m=0}^{\infty}\frac{(-1)^{m}(2m+1)e^{-\pi t(2m+1)}G_{r-N,t,k}}{t^2(2m+1)^2+(r+c)^2} &\\&+\frac{(-1)^{c-1}t\sinh(\pi t)}{2\sqrt{k}\pi}\sum\limits_{r=1}^{\infty} \sum\limits_{m=0}^{\infty}\frac{(-1)^{m}(2m+1)e^{-\pi t(2m+1)}G_{-r-N,t,k}}{t^2(2m+1)^2+(r-c)^2} . & 
        \label{q(a+b)formula} \tag{4.7}
    \end{flalign*}

    And the counterpart for $N+c \leq 0$ is

    \begin{flalign*}
         0 & = U_{N+c}- \frac{(N+c)\pi^2\coth(\pi t)}{3k} \sum\limits_{r=1}^{\infty} \frac{(-1)^{r-1}e^{-2\pi t r}}{4t^2r^2+(N+c)^2} &\\& - t\coth(\pi t) \sum\limits_{r=1}^{\infty} \frac{(-1)^{r-1}re^{-2\pi t r}}{4t^2r^2+(N+c)^2} - \frac{\coth(\pi t)}{\pi}\sum\limits_{r=1}^{\infty} \frac{(-1)^{r-1}e^{-2\pi t r}(4t^2r^2-(N+c)^2)}{(4t^2r^2+(N+c)^2)^2} + \frac{\sinh(\pi t) G_{-N,t,k}}{8\sqrt{k} t \cosh(\frac{\pi c}{2t})} &\\& 
        + \frac{(-1)^{c-1}t\sinh(\pi t)G_{-N,t,k}}{2\sqrt{k}\pi}\sum\limits_{r=0}^{\infty}\frac{(-1)^{r}(2r+1)e^{-\pi t(2r+1)}}{t^2(2r+1)^2+c^2}
        + \frac{\sinh(\pi t)}{8\sqrt{k}t}\sum\limits_{r=1}^{\infty}(-1)^r\frac{G_{r-N,t,k}}{ \cosh(\frac{\pi (r+c)}{2t})} &\\&+\frac{\sinh(\pi t)}{8\sqrt{k}t}\sum\limits_{r=1}^{\infty}(-1)^r\frac{G_{-r-N,t,k}}{ \cosh(\frac{\pi (r-c)}{2t})}  +\frac{(-1)^{c-1}t\sinh(\pi t)}{2\sqrt{k} \pi}\sum\limits_{r=1}^{\infty} \sum\limits_{m=0}^{\infty}\frac{(-1)^{m}(2m+1)e^{-\pi t(2m+1)}G_{r-N,t,k}}{t^2(2m+1)^2+(r+c)^2} &\\&+\frac{(-1)^{c-1}t\sinh(\pi t)}{2\sqrt{k}\pi}\sum\limits_{r=1}^{\infty} \sum\limits_{m=0}^{\infty}\frac{(-1)^{m}(2m+1)e^{-\pi t(2m+1)}G_{-r-N,t,k}}{t^2(2m+1)^2+(r-c)^2} ,
        \label{0(a+b)formula} \tag{4.8}
    \end{flalign*}

    where 
    \[
    U_R = \begin{cases}
        \frac{ \pi^2}{12} - \frac{\pi \coth(\pi t)}{4} + \frac{\pi \coth(\pi t)}{48 t^2} & \textit{, if } R = 0, \vspace{0.7cm}\\
        \frac{\pi^2}{3Rk(e^{2\pi t}-1)}+ \frac{\pi - \coth(\pi t)}{2\pi R^2}+ \frac{(-1)^R(1-t^2)\coth(\pi t)}{\pi R^2} + \frac{(-1)^{R-1} \pi^3 \coth(\pi t)}{12kt \sinh(\frac{\pi R}{2t})} \\+ \frac{(-1)^R\coth(\pi t)(t-\frac{1}{t})}{4R\sinh(\frac{\pi R}{2t})}+ \frac{(-1)^R\pi \cosh(\frac{\pi R}{2t})\coth(\pi t)}{8\sinh^2(\frac{\pi R}{2t})} & \textit{, if else.}
    \end{cases}
    \]

    With all the equalities being valid for $t>0$.
\vspace{0.7cm}
\section{Analytic expressions for certain general arithmetical sums}

Let $h$ be a function such that for all natural numbers $m$, we have $|h(m)|\leq M$ for some real value $M$, and in addition the series $\sum\limits_{r=1}^{\infty} \frac{h(r)}{r+z}$ converges absolutely for all $z$, where the sum excludes the term $r=-z$, when $z$ is a negative integer. Then we propose to study the sums

\begin{flalign*}
    \sum\limits_{a=1}^{N-1} \frac{h(a) q_k(N-a)}{(N-a)^2}, \label{finitesum} \tag{5.1}\\
    \sum\limits_{a=1}^{\infty} \frac{h(a) q_k(N+a)}{(N+a)^2}.\label{infinitesum} \tag{5.2}
\end{flalign*}

First we shall establish the following lemma, which will then allow us to express the above sums in an analytic fashion.
\vspace{0.7cm}
\begin{lemma} For a function $h$ satisfying the conditions described, The following double and triple sums are absolutely convergent for real $t > 0$.

\begin{flalign*}
    &\sum\limits_{r=1}^{\infty} \sum\limits_{a=1}^{\infty}\frac{(-1)^r h(a)G_{r-N,t,k}}{ \cosh(\frac{\pi (r-a)}{2t})},  \hspace{0.3cm}
    \sum\limits_{r=1}^{\infty} \sum\limits_{a=1}^{\infty}\frac{(-1)^r h(a)G_{r-N,t,k}}{ \cosh(\frac{\pi (r+a)}{2t})}, \\
    &\sum\limits_{r=1}^{\infty} \sum\limits_{a=1}^{\infty}\frac{(-1)^r h(a)G_{-r-N,t,k}}{ \cosh(\frac{\pi (r+a)}{2t})},  \hspace{0.3cm}
    \sum\limits_{r=1}^{\infty} \sum\limits_{a=1}^{\infty}\frac{(-1)^r h(a)G_{-r-N,t,k}}{ \cosh(\frac{\pi (r-a)}{2t})}, \\
    &\sum\limits_{r=1}^{\infty} \sum\limits_{a=1}^{\infty} \sum\limits_{m=0}^{\infty}\frac{(-1)^{m}(2m+1)h(a)e^{-\pi t(2m+1)}G_{r-N,t,k}}{t^2(2m+1)^2+(r-a)^2},  \hspace{0.3cm} 
    \sum\limits_{r=1}^{\infty} \sum\limits_{a=1}^{\infty} \sum\limits_{m=0}^{\infty}\frac{(-1)^{m}(2m+1)h(a)e^{-\pi t(2m+1)}G_{r-N,t,k}}{t^2(2m+1)^2+(r+a)^2}, \\
    &\sum\limits_{r=1}^{\infty} \sum\limits_{a=1}^{\infty} \sum\limits_{m=0}^{\infty}\frac{(-1)^{m}(2m+1)h(a)e^{-\pi t(2m+1)}G_{-r-N,t,k}}{t^2(2m+1)^2+(r+a)^2},  \hspace{0.3cm} 
    \sum\limits_{r=1}^{\infty} \sum\limits_{a=1}^{\infty} \sum\limits_{m=0}^{\infty}\frac{(-1)^{m}(2m+1)h(a)e^{-\pi t(2m+1)}G_{-r-N,t,k}}{t^2(2m+1)^2+(r-a)^2}.
\end{flalign*}
    
\end{lemma}

\begin{proof}[Proof of Lemma 13]
    \begin{flalign*}
        \sum\limits_{r=1}^{\infty} \sum\limits_{a=1}^{\infty} \left| \frac{(-1)^r h(a)G_{r-N,t,k}}{ \cosh(\frac{\pi (r-a)}{2t})} \right| &\leq M \sum\limits_{r=1}^{\infty} \sum\limits_{a=1}^{\infty}  \frac{|G_{r-N,t,k}|}{ \cosh(\frac{\pi (r-a)}{2t})} \\ 
        & = M \sum\limits_{r=1}^{\infty}|G_{r-N,t,k}| \sum\limits_{a=1}^{\infty}  \frac{1}{ \cosh(\frac{\pi (r-a)}{2t})} \\ 
        & = M \sum\limits_{r=1}^{\infty}|G_{r-N,t,k}| \sum\limits_{a=0}^{r-1}  \frac{1}{ \cosh(\frac{\pi a}{2t})} +
        M \sum\limits_{r=1}^{\infty}|G_{r-N,t,k}| \sum\limits_{a=1}^{\infty}  \frac{1}{ \cosh(\frac{\pi a}{2t})}, \label{proofeq34} \tag{5.3}
    \end{flalign*}
    
    \begin{flalign*}
        \sum\limits_{r=1}^{\infty} \sum\limits_{a=1}^{\infty} \left| \frac{(-1)^r h(a)G_{-r-N,t,k}}{ \cosh(\frac{\pi (r-a)}{2t})} \right| &\leq M \sum\limits_{r=1}^{\infty} \sum\limits_{a=1}^{\infty}  \frac{|G_{-r-N,t,k}|}{ \cosh(\frac{\pi (r-a)}{2t})} \\ 
        & = M \sum\limits_{r=1}^{\infty}|G_{-r-N,t,k}| \sum\limits_{a=1}^{\infty}  \frac{1}{ \cosh(\frac{\pi (r-a)}{2t})} \\ 
        & = M \sum\limits_{r=1}^{\infty}|G_{-r-N,t,k}| \sum\limits_{a=0}^{r-1}  \frac{1}{ \cosh(\frac{\pi a}{2t})} +
        M \sum\limits_{r=1}^{\infty}|G_{-r-N,t,k}| \sum\limits_{a=1}^{\infty}  \frac{1}{ \cosh(\frac{\pi a}{2t})}. \label{proofeq35} \tag{5.4}
    \end{flalign*}
    Since $\sum\limits_{a=0}^{r-1}  \frac{1}{ \cosh(\frac{\pi a}{2t})}$ is bounded for all natural numbers $r$, and since by the proof of \textbf{Lemma 7}, $\sum\limits_{r=1}^{\infty}|G_{r-N,t,k}|$ and $\sum\limits_{r=1}^{\infty}|G_{-r-N,t,k}|$ converges, we must have $\sum\limits_{r=1}^{\infty} \sum\limits_{a=1}^{\infty}\frac{(-1)^r h(a)G_{r-N,t,k}}{ \cosh(\frac{\pi (r-a)}{2t})}$ and $\sum\limits_{r=1}^{\infty} \sum\limits_{a=1}^{\infty}\frac{(-1)^r h(a)G_{-r-N,t,k}}{ \cosh(\frac{\pi (r-a)}{2t})}$ absolutely convergent. Moving on, we have 
    \begin{flalign*}
        \sum\limits_{r=1}^{\infty} \sum\limits_{a=1}^{\infty}\left|\frac{(-1)^r h(a)G_{r-N,t,k}}{ \cosh(\frac{\pi (r+a)}{2t})} \right| & \leq M\sum\limits_{r=1}^{\infty} \sum\limits_{a=1}^{\infty}\frac{|G_{r-N,t,k}|}{ \cosh(\frac{\pi (r+a)}{2t})} \\ &  <   2M\sum\limits_{r=1}^{\infty} e^{-\frac{\pi r}{2t}}|G_{r-N,t,k}|\sum\limits_{a=1}^{\infty} e^{-\frac{\pi a}{2t}} , \label{proofeq36} \tag{5.5}
    \end{flalign*}

    \begin{flalign*}
        \sum\limits_{r=1}^{\infty} \sum\limits_{a=1}^{\infty}\left|\frac{(-1)^r h(a)G_{-r-N,t,k}}{ \cosh(\frac{\pi (r+a)}{2t})} \right| & \leq M\sum\limits_{r=1}^{\infty} \sum\limits_{a=1}^{\infty}\frac{|G_{-r-N,t,k}|}{ \cosh(\frac{\pi (r+a)}{2t})} \\ &  <   2M\sum\limits_{r=1}^{\infty} e^{-\frac{\pi r}{2t}}|G_{-r-N,t,k}|\sum\limits_{a=1}^{\infty} e^{-\frac{\pi a}{2t}} .
        \label{proofeq37} \tag{5.6}
    \end{flalign*}
    Thus we have $\sum\limits_{r=1}^{\infty} \sum\limits_{a=1}^{\infty}\frac{(-1)^r h(a)G_{r-N,t,k}}{ \cosh(\frac{\pi (r+a)}{2t})}$ and $\sum\limits_{r=1}^{\infty} \sum\limits_{a=1}^{\infty}\frac{(-1)^r h(a)G_{-r-N,t,k}}{ \cosh(\frac{\pi (r+a)}{2t})}$ absolutely convergent. Next we have

    \begin{flalign*}
        \sum\limits_{r=1}^{\infty} \sum\limits_{a=1}^{\infty} \sum\limits_{m=0}^{\infty}\left|\frac{(-1)^{m}(2m+1)h(a)e^{-\pi t(2m+1)}G_{r-N,t,k}}{t^2(2m+1)^2+(r-a)^2} \right| &\leq M\sum\limits_{r=1}^{\infty} \sum\limits_{a=1}^{\infty} \sum\limits_{m=0}^{\infty}\frac{(2m+1)e^{-\pi t(2m+1)}|G_{r-N,t,k}|}{t^2(2m+1)^2+(r-a)^2} \\ & \leq 
        M\sum\limits_{r=1}^{\infty} \sum\limits_{a=1}^{\infty} \frac{|G_{r-N,t,k}|}{t^2+(r-a)^2} \sum\limits_{m=0}^{\infty}(2m+1)e^{-\pi t(2m+1)} \\ & = 
        \left(M\sum\limits_{m=0}^{\infty}(2m+1)e^{-\pi t(2m+1)}\right)\sum\limits_{r=1}^{\infty}|G_{r-N,t,k}| \sum\limits_{a=0}^{r-1} \frac{1}{t^2+a^2} 
        \\& + \left(M\sum\limits_{m=0}^{\infty}(2m+1)e^{-\pi t(2m+1)}\right)\sum\limits_{r=1}^{\infty}|G_{r-N,t,k}| \sum\limits_{a=1}^{\infty} \frac{1}{t^2+a^2}, \label{proofeq38} \tag{5.7}
    \end{flalign*}

    \begin{flalign*}
        \sum\limits_{r=1}^{\infty} \sum\limits_{a=1}^{\infty} \sum\limits_{m=0}^{\infty}\left|\frac{(-1)^{m}(2m+1)h(a)e^{-\pi t(2m+1)}G_{-r-N,t,k}}{t^2(2m+1)^2+(r-a)^2} \right| &\leq M\sum\limits_{r=1}^{\infty} \sum\limits_{a=1}^{\infty} \sum\limits_{m=0}^{\infty}\frac{(2m+1)e^{-\pi t(2m+1)}|G_{-r-N,t,k}|}{t^2(2m+1)^2+(r-a)^2} \\ & \leq 
        M\sum\limits_{r=1}^{\infty} \sum\limits_{a=1}^{\infty} \frac{|G_{-r-N,t,k}|}{t^2+(r-a)^2} \sum\limits_{m=0}^{\infty}(2m+1)e^{-\pi t(2m+1)} \\ & = 
        \left(M\sum\limits_{m=0}^{\infty}(2m+1)e^{-\pi t(2m+1)}\right)\sum\limits_{r=1}^{\infty}|G_{-r-N,t,k}| \sum\limits_{a=0}^{r-1} \frac{1}{t^2+a^2} 
        \\& + \left(M\sum\limits_{m=0}^{\infty}(2m+1)e^{-\pi t(2m+1)}\right)\sum\limits_{r=1}^{\infty}|G_{-r-N,t,k}| \sum\limits_{a=1}^{\infty} \frac{1}{t^2+a^2}. \label{proofeq39} \tag{5.8}
    \end{flalign*}

    Since $\sum\limits_{a=0}^{r-1} \frac{1}{t^2+a^2}$ is bounded for all natural numbers $r$, and since $\sum\limits_{r=1}^{\infty}|G_{r-N,t,k}|$ and $\sum\limits_{r=1}^{\infty}|G_{-r-N,t,k}|$ converges, we have $\sum\limits_{r=1}^{\infty} \sum\limits_{a=1}^{\infty} \sum\limits_{m=0}^{\infty}\frac{(-1)^{m}(2m+1)h(a)e^{-\pi t(2m+1)}G_{r-N,t,k}}{t^2(2m+1)^2+(r-a)^2}$ and $ 
    \sum\limits_{r=1}^{\infty} \sum\limits_{a=1}^{\infty} \sum\limits_{m=0}^{\infty}\frac{(-1)^{m}(2m+1)h(a)e^{-\pi t(2m+1)}G_{-r-N,t,k}}{t^2(2m+1)^2+(r-a)^2}$ absolutely convergent. Lastly we move on to

    \begin{flalign*}
        &\sum\limits_{r=1}^{\infty} \sum\limits_{a=1}^{\infty} \sum\limits_{m=0}^{\infty}\left|\frac{(-1)^{m}(2m+1)h(a)e^{-\pi t(2m+1)}G_{r-N,t,k}}{t^2(2m+1)^2+(r+a)^2} \right| \\&\leq M\sum\limits_{r=1}^{\infty} \sum\limits_{a=1}^{\infty} \sum\limits_{m=0}^{\infty} \frac{(2m+1)e^{-\pi t(2m+1)}|G_{r-N,t,k}|}{t^2(2m+1)^2+(r+a)^2} \\ &\leq M \left( \sum\limits_{m=0}^{\infty} (2m+1)e^{-\pi t(2m+1)}\right)\left( \sum\limits_{r=1}^{\infty} \sum\limits_{a=1}^{\infty}  \frac{|G_{r-N,t,k}|}{t^2+(r+a)^2}\right) \\ & < M \left( \sum\limits_{m=0}^{\infty} (2m+1)e^{-\pi t(2m+1)}\right)
        \left( \sum\limits_{a=1}^{\infty} \frac{1}{t^2+a^2} \right)\left( \sum\limits_{r=1}^{\infty}  |G_{r-N,t,k}|\right), 
        \label{proofeq40} \tag{5.9}
    \end{flalign*}

    \begin{flalign*}
        &\sum\limits_{r=1}^{\infty} \sum\limits_{a=1}^{\infty} \sum\limits_{m=0}^{\infty}\left|\frac{(-1)^{m}(2m+1)h(a)e^{-\pi t(2m+1)}G_{-r-N,t,k}}{t^2(2m+1)^2+(r+a)^2} \right| \\&\leq M\sum\limits_{r=1}^{\infty} \sum\limits_{a=1}^{\infty} \sum\limits_{m=0}^{\infty} \frac{(2m+1)e^{-\pi t(2m+1)}|G_{-r-N,t,k}|}{t^2(2m+1)^2+(r+a)^2} \\ &\leq M \left( \sum\limits_{m=0}^{\infty} (2m+1)e^{-\pi t(2m+1)}\right)\left( \sum\limits_{r=1}^{\infty} \sum\limits_{a=1}^{\infty}  \frac{|G_{-r-N,t,k}|}{t^2+(r+a)^2}\right) \\ & < M \left( \sum\limits_{m=0}^{\infty} (2m+1)e^{-\pi t(2m+1)}\right)
        \left( \sum\limits_{a=1}^{\infty} \frac{1}{t^2+a^2} \right)\left( \sum\limits_{r=1}^{\infty}  |G_{-r-N,t,k}|\right).
        \label{proofeq41} \tag{5.10}
    \end{flalign*}

    Therefore we have $\sum\limits_{r=1}^{\infty} \sum\limits_{a=1}^{\infty} \sum\limits_{m=0}^{\infty}\frac{(-1)^{m}(2m+1)h(a)e^{-\pi t(2m+1)}G_{r-N,t,k}}{t^2(2m+1)^2+(r+a)^2}$ and $\sum\limits_{r=1}^{\infty} \sum\limits_{a=1}^{\infty} \sum\limits_{m=0}^{\infty}\frac{(-1)^{m}(2m+1)h(a)e^{-\pi t(2m+1)}G_{-r-N,t,k}}{t^2(2m+1)^2+(r+a)^2}$ absolutely convergent.

\end{proof}

Next, we proceed to formulate an analytic representation for the sums \eqref{finitesum} and \eqref{infinitesum}. 
\vspace{0.7cm}
\begin{proposition} For a function $h$ satisfying the conditions described, the following holds for all natural numbers $N$ and real $t > 0$.

\begin{flalign*}
    &\sum\limits_{a=1}^{N-1} \frac{h(a) q_k(N-a)}{(N-a)^2} 
    \\ &= \left(\frac{ \pi^2}{12} - \frac{\pi \coth(\pi t)}{4} + \frac{\pi \coth(\pi t)}{48 t^2} \right) h(N)+ \frac{\pi^2}{3k(e^{2\pi t}-1)}{\sum\limits_{a=1}^{\infty} }^*\frac{h(a)}{N-a}+ \frac{\pi - \coth(\pi t)}{2\pi}{\sum\limits_{a=1}^{\infty} }^*\frac{h(a)}{(N-a)^2} &\\&+ \frac{(-1)^{N}(1-t^2)\coth(\pi t)}{\pi}{\sum\limits_{a=1}^{\infty} }^* \frac{(-1)^ah(a)}{(N-a)^2} + \frac{(-1)^{N} \pi^3 \coth(\pi t)}{12kt} {\sum\limits_{a=1}^{\infty}}^* \frac{(-1)^{a-1}h(a)}{\sinh(\frac{\pi (N-a)}{2t})}  &\\&+ \frac{(-1)^{N}\coth(\pi t)(t-\frac{1}{t})}{4} {\sum\limits_{a=1}^{\infty}}^* \frac{(-1)^ah(a)}{(N-a)\sinh(\frac{\pi (N-a)}{2t})} + \frac{(-1)^{N}\pi \coth(\pi t)}{8} {\sum\limits_{a=1}^{\infty}}^* \frac{(-1)^ah(a)\cosh(\frac{\pi(N-a)}{2t})}{\sinh^2(\frac{\pi (N-a)}{2t})} &\\&- \frac{\pi^2\coth(\pi t)}{3k} \sum\limits_{r=1}^{\infty} \sum\limits_{a=1}^{\infty} \frac{(-1)^{r-1}e^{-2\pi t r}(N-a)h(a)}{4t^2r^2+(N-a)^2}  - t\coth(\pi t) \sum\limits_{r=1}^{\infty} \sum\limits_{a=1}^{\infty} \frac{(-1)^{r-1}re^{-2\pi t r}h(a)}{4t^2r^2+(N-a)^2} &\\& - \frac{\coth(\pi t)}{\pi}\sum\limits_{r=1}^{\infty}\sum\limits_{a=1}^{\infty} \frac{(-1)^{r-1}e^{-2\pi t r}(4t^2r^2-(N-a)^2)h(a)}{(4t^2r^2+(N-a)^2)^2} + \frac{\sinh(\pi t) G_{-N,t,k}}{8\sqrt{k} t }\sum\limits_{a=1}^{\infty} \frac{h(a)}{\cosh(\frac{\pi a}{2t})} &\\& + \frac{t\sinh(\pi t)G_{-N,t,k}}{2\sqrt{k}\pi}\sum\limits_{r=0}^{\infty} \sum\limits_{a=1}^{\infty}\frac{(-1)^{r+a-1}(2r+1)e^{-\pi t(2r+1)}h(a)}{t^2(2r+1)^2+a^2} &\\& +\frac{\sinh(\pi t)}{8\sqrt{k}t}\sum\limits_{r=1}^{\infty}\sum\limits_{a=1}^{\infty}\frac{(-1)^rG_{r-N,t,k}h(a)}{ \cosh(\frac{\pi (r-a)}{2t})}+\frac{\sinh(\pi t)}{8\sqrt{k}t}\sum\limits_{r=1}^{\infty}\sum\limits_{a=1}^{\infty}\frac{(-1)^rG_{-r-N,t,k}h(a)}{ \cosh(\frac{\pi (r+a)}{2t})} &\\& +\frac{t\sinh(\pi t)}{2\sqrt{k} \pi}\sum\limits_{r=1}^{\infty} \sum\limits_{m=0}^{\infty}\sum\limits_{a=1}^{\infty} \frac{(-1)^{m+a-1}(2m+1)e^{-\pi t(2m+1)}G_{r-N,t,k}h(a)}{t^2(2m+1)^2+(r-a)^2} &\\&+\frac{t\sinh(\pi t)}{2\sqrt{k}\pi}\sum\limits_{r=1}^{\infty} \sum\limits_{m=0}^{\infty}\sum\limits_{a=1}^{\infty}\frac{(-1)^{m+a-1}(2m+1)e^{-\pi t(2m+1)}G_{-r-N,t,k}h(a)}{t^2(2m+1)^2+(r+a)^2} ,& \label{analyticStatement1} \tag{5.11}
\end{flalign*}
where $*$ indicates that the term inside the summation exhibiting  a singularity is taken as zero, and
\begin{flalign*}
    &\sum\limits_{a=1}^{\infty} \frac{h(a) q_k(N+a)}{(N+a)^2} &\\& =  \frac{\pi^2}{3k(e^{2\pi t}-1)} \sum\limits_{a=1}^{\infty} \frac{h(a)}{N+a}+ \frac{\pi - \coth(\pi t)}{2\pi}\sum\limits_{a=1}^{\infty} \frac{h(a)}{(N+a)^2}+ \frac{(-1)^{N}(1-t^2)\coth(\pi t)}{\pi }\sum\limits_{a=1}^{\infty} \frac{(-1)^ah(a)}{(N+a)^2} &\\& + \frac{(-1)^{N} \pi^3 \coth(\pi t)}{12kt }\sum\limits_{a=1}^{\infty} \frac{(-1)^{a-1}h(a)}{\sinh(\frac{\pi (N+a)}{2t})}  + \frac{(-1)^{N}\coth(\pi t)(t-\frac{1}{t})}{4}\sum\limits_{a=1}^{\infty} \frac{(-1)^ah(a)}{(N+a)\sinh(\frac{\pi (N+a)}{2t})}&\\&+ \frac{(-1)^{N}\pi \coth(\pi t)}{8}\sum\limits_{a=1}^{\infty} \frac{(-1)^ah(a)\cosh(\frac{\pi(N+a)}{2t})}{\sinh^2(\frac{\pi (N+a)}{2t})} - \frac{\pi^2\coth(\pi t)}{3k} \sum\limits_{r=1}^{\infty}\sum\limits_{a=1}^{\infty} \frac{(-1)^{r-1}e^{-2\pi t r}(N+a)h(a)}{4t^2r^2+(N+a)^2} &\\& - t\coth(\pi t) \sum\limits_{r=1}^{\infty} \sum\limits_{a=1}^{\infty} \frac{(-1)^{r-1}re^{-2\pi t r}h(a)}{4t^2r^2+(N+a)^2} - \frac{\coth(\pi t)}{\pi}\sum\limits_{r=1}^{\infty} \sum\limits_{a=1}^{\infty} \frac{(-1)^{r-1}e^{-2\pi t r}(4t^2r^2-(N+a)^2)h(a)}{(4t^2r^2+(N+a)^2)^2} &\\& + \frac{\sinh(\pi t) G_{-N,t,k}}{8\sqrt{k} t } \sum\limits_{a=1}^{\infty} \frac{h(a)}{\cosh(\frac{\pi a}{2t})} +\frac{t\sinh(\pi t)G_{-N,t,k}}{2\sqrt{k}\pi}\sum\limits_{r=0}^{\infty} \sum\limits_{a=1}^{\infty}\frac{(-1)^{r+a-1}(2r+1)e^{-\pi t(2r+1)}h(a)}{t^2(2r+1)^2+a^2} &\\&+\frac{\sinh(\pi t)}{8\sqrt{k}t}\sum\limits_{r=1}^{\infty} \sum\limits_{a=1}^{\infty}\frac{(-1)^rG_{r-N,t,k}h(a)}{ \cosh(\frac{\pi (r+a)}{2t})} +\frac{\sinh(\pi t)}{8\sqrt{k}t}\sum\limits_{r=1}^{\infty} \sum\limits_{a=1}^{\infty}(-1)^r\frac{G_{-r-N,t,k}h(a)}{ \cosh(\frac{\pi (r-a)}{2t})} &\end{flalign*}\begin{flalign*}& +\frac{t\sinh(\pi t)}{2\sqrt{k} \pi}\sum\limits_{r=1}^{\infty} \sum\limits_{m=0}^{\infty}\sum\limits_{a=1}^{\infty}\frac{(-1)^{m+a-1}(2m+1)e^{-\pi t(2m+1)}G_{r-N,t,k}h(a)}{t^2(2m+1)^2+(r+a)^2} &\\&+\frac{t\sinh(\pi t)}{2\sqrt{k}\pi}\sum\limits_{r=1}^{\infty} \sum\limits_{m=0}^{\infty}\sum\limits_{a=1}^{\infty}\frac{(-1)^{m+a-1}(2m+1)e^{-\pi t(2m+1)}G_{-r-N,t,k}h(a)}{t^2(2m+1)^2+(r-a)^2} .& \label{analyticStatement2} \tag{5.12}
\end{flalign*}
\end{proposition}

\begin{proof}[Proof of Proposition 14]
    For the first equality \eqref{analyticStatement1}, we note that $\sum\limits_{a=1}^{N-1} \frac{h(a) q_k(N-a)}{(N-a)^2}= {\sum\limits_{a=1}^{\infty} }^* \frac{h(a) q_k(N-a)}{(N-a)^2}$, (whereas the value of $\frac{q_k(N)}{N^2}$ is to be regarded as zero for $N\leq 0$) and substitute the analytic expression for $q_k(N-a)$. The second equality \eqref{analyticStatement2} is arrived in the same manner substituting for $q_k(N+a)$. The interchanging of the order of summation in double and triple sums is justified by \textbf{Lemma 13}.
\end{proof}

As an application of the above formula, we present the analytic representation for a general arithmetical sum whose terms are summed over all the positive integer solutions $(a,b)$ of the Diophantine equation $N= da^2+kb^2$, where $d, k \in \mathbb{N}$ . This is achieved when we put $h(m) = q_d(m)g(\sqrt{\frac{m}{d}}) $ in the equalities \eqref{analyticStatement1} and \eqref{analyticStatement2}, where $g$ is some function such that for all natural numbers $m$, we have $|g(m)| \leq M$ for some real value $M$ and in addition the series ${\sum\limits_{r=1}^{\infty}}^* \frac{g(r)}{z+r^2}$ converges absolutely for all $z$. Then the sum $\sum\limits_{a=1}^{N-1} \frac{h(a) q_k(N-a)}{(N-a)^2}$ becomes $\sum\limits_{a=1}^{N-1} \frac{g(\sqrt{\frac{a}{d}}) q_d(a) q_k(N-a)}{(N-a)^2} = \frac{1}{k^2}\sum\limits_{\substack{da^2+kb^2=N \\ (a,b) \in \mathbb{N}^2}} \frac{g(a)}{b^4}$. Therefore we state
\vspace{0.7cm}
\begin{proposition} Let $g$ be a function such that for all natural numbers $m$, we have $|g(m)| \leq M$ for some real value $M$. And in addition the series ${\sum\limits_{r=1}^{\infty}}^* \frac{g(r)}{z+r^2}$ converges absolutely for all $z$. Then for all natural numbers $N$ and real $t>0$, we have
    \begin{flalign*}
    &\frac{1}{k^2}\sum\limits_{\substack{da^2+kb^2=N \\ (a,b) \in \mathbb{N}^2}} \frac{g(a)}{b^4} 
    \\ &= \left(\frac{ \pi^2}{12} - \frac{\pi \coth(\pi t)}{4} + \frac{\pi \coth(\pi t)}{48 t^2} \right)\epsilon_N + \frac{\pi^2}{3k(e^{2\pi t}-1)}{\sum\limits_{a=1}^{\infty} }^*\frac{g(a)}{N-da^2}+ \frac{\pi - \coth(\pi t)}{2\pi}{\sum\limits_{a=1}^{\infty} }^*\frac{g(a)}{(N-da^2)^2} &\\&+ \frac{(-1)^{N}(1-t^2)\coth(\pi t)}{\pi}{\sum\limits_{a=1}^{\infty} }^* \frac{(-1)^{ad}g(a)}{(N-da^2)^2} + \frac{(-1)^{N} \pi^3 \coth(\pi t)}{12kt} {\sum\limits_{a=1}^{\infty}}^* \frac{(-1)^{ad-1}g(a)}{\sinh(\frac{\pi (N-da^2)}{2t})}  &\end{flalign*} \begin{flalign*}&+ \frac{(-1)^{N}\coth(\pi t)(t-\frac{1}{t})}{4} {\sum\limits_{a=1}^{\infty}}^* \frac{(-1)^{ad}g(a)}{(N-da^2)\sinh(\frac{\pi (N-da^2)}{2t})} + \frac{(-1)^{N}\pi \coth(\pi t)}{8} {\sum\limits_{a=1}^{\infty}}^* \frac{(-1)^{ad}g(a)\cosh(\frac{\pi(N-da^2)}{2t})}{\sinh^2(\frac{\pi (N-da^2)}{2t})} &\\&- \frac{\pi^2\coth(\pi t)}{3k} \sum\limits_{r=1}^{\infty} \sum\limits_{a=1}^{\infty} \frac{(-1)^{r-1}e^{-2\pi t r}(N-da^2)g(a)}{4t^2r^2+(N-da^2)^2}  - t\coth(\pi t) \sum\limits_{r=1}^{\infty} \sum\limits_{a=1}^{\infty} \frac{(-1)^{r-1}re^{-2\pi t r}g(a)}{4t^2r^2+(N-da^2)^2} &\\& - \frac{\coth(\pi t)}{\pi}\sum\limits_{r=1}^{\infty}\sum\limits_{a=1}^{\infty} \frac{(-1)^{r-1}e^{-2\pi t r}(4t^2r^2-(N-da^2)^2)g(a)}{(4t^2r^2+(N-da^2)^2)^2} + \frac{\sinh(\pi t) G_{-N,t,k}}{8\sqrt{k} t }\sum\limits_{a=1}^{\infty} \frac{g(a)}{\cosh(\frac{\pi da^2}{2t})} &\\& + \frac{t\sinh(\pi t)G_{-N,t,k}}{2\sqrt{k}\pi}\sum\limits_{r=0}^{\infty} \sum\limits_{a=1}^{\infty}\frac{(-1)^{r+da-1}(2r+1)e^{-\pi t(2r+1)}g(a)}{t^2(2r+1)^2+d^2a^4}+\frac{\sinh(\pi t)}{8\sqrt{k}t}\sum\limits_{r=1}^{\infty}\sum\limits_{a=1}^{\infty}\frac{(-1)^rG_{r-N,t,k}g(a)}{ \cosh(\frac{\pi (r-da^2)}{2t})}&\\&+\frac{\sinh(\pi t)}{8\sqrt{k}t}\sum\limits_{r=1}^{\infty}\sum\limits_{a=1}^{\infty}\frac{(-1)^rG_{-r-N,t,k}g(a)}{ \cosh(\frac{\pi (r+da^2)}{2t})}  +\frac{t\sinh(\pi t)}{2\sqrt{k} \pi}\sum\limits_{r=1}^{\infty} \sum\limits_{m=0}^{\infty}\sum\limits_{a=1}^{\infty} \frac{(-1)^{m+da-1}(2m+1)e^{-\pi t(2m+1)}G_{r-N,t,k}g(a)}{t^2(2m+1)^2+(r-da^2)^2} &\\&+\frac{t\sinh(\pi t)}{2\sqrt{k}\pi}\sum\limits_{r=1}^{\infty} \sum\limits_{m=0}^{\infty}\sum\limits_{a=1}^{\infty}\frac{(-1)^{m+da-1}(2m+1)e^{-\pi t(2m+1)}G_{-r-N,t,k}g(a)}{t^2(2m+1)^2+(r+da^2)^2} , \label{analyticStatementForSumofSquares} \tag{5.13}
\end{flalign*}
where $d,k \in \mathbb{N}$ and 
\[
\epsilon_N = \begin{cases}
    g(m)&, \textit{ if } N = dm^2 \textit{ for some natural number }m, \\
    0&, \textit{ if else.}
\end{cases}
\]
\end{proposition}

Next, we state the counterpart of the formula \eqref{analyticStatementForSumofSquares} for a general arithmetical sum whose terms are summed over all the positive integer solutions of $(a,b)$ of the Diophantine equation $kb^2-da^2=N$, where $d, k \in \mathbf{N}$. Again we put $h(m) = q_d(m)g(\sqrt{\frac{m}{d}}) $ in eq, where $g$ is some function such that for all natural numbers $m$, we have $|g(m)| \leq M$ for some real value $M$ and in addition the series ${\sum\limits_{r=1}^{\infty}}^* \frac{g(r)}{z+r^2}$ converges absolutely for all $z$. Then the sum $\sum\limits_{a=1}^{\infty} \frac{h(a) q_k(N+a)}{(N+a)^2}$ becomes $\sum\limits_{a=1}^{\infty} \frac{g(\sqrt{\frac{a}{d}}) q_d(a) q_k(N+a)}{(N+a)^2}= \frac{1}{k^2}\sum\limits_{\substack{kb^2-da^2=N \\ (a,b) \in \mathbb{N}^2}} \frac{g(a)}{b^4}$. Thus, we state
\vspace{0.7cm}
\begin{proposition}Let $g$ be a function such that for all natural numbers $m$, we have $|g(m)| \leq M$ for some real value $M$. And in addition the series ${\sum\limits_{r=1}^{\infty}}^* \frac{g(r)}{z+r^2}$ converges absolutely for all $z$. Then for all natural numbers $N$ and real $t>0$, we have

\begin{flalign*}
    &\frac{1}{k^2}\sum\limits_{\substack{kb^2-da^2=N \\ (a,b) \in \mathbb{N}^2}} \frac{g(a)}{b^4}  &\\& =  \frac{\pi^2}{3k(e^{2\pi t}-1)} \sum\limits_{a=1}^{\infty} \frac{g(a)}{N+da^2}+ \frac{\pi - \coth(\pi t)}{2\pi}\sum\limits_{a=1}^{\infty} \frac{g(a)}{(N+da^2)^2}+ \frac{(-1)^{N}(1-t^2)\coth(\pi t)}{\pi }\sum\limits_{a=1}^{\infty} \frac{(-1)^{ad}g(a)}{(N+da^2)^2} &\\& + \frac{(-1)^{N} \pi^3 \coth(\pi t)}{12kt }\sum\limits_{a=1}^{\infty} \frac{(-1)^{ad-1}g(a)}{\sinh(\frac{\pi (N+da^2)}{2t})}  + \frac{(-1)^{N}\coth(\pi t)(t-\frac{1}{t})}{4}\sum\limits_{a=1}^{\infty} \frac{(-1)^{ad}g(a)}{(N+da^2)\sinh(\frac{\pi (N+da^2)}{2t})}&\\&+ \frac{(-1)^{N}\pi \coth(\pi t)}{8}\sum\limits_{a=1}^{\infty} \frac{(-1)^{ad}g(a)\cosh(\frac{\pi(N+da^2)}{2t})}{\sinh^2(\frac{\pi (N+da^2)}{2t})} - \frac{\pi^2\coth(\pi t)}{3k} \sum\limits_{r=1}^{\infty}\sum\limits_{a=1}^{\infty} \frac{(-1)^{r-1}e^{-2\pi t r}(N+da^2)g(a)}{4t^2r^2+(N+da^2)^2} &\\& - t\coth(\pi t) \sum\limits_{r=1}^{\infty} \sum\limits_{a=1}^{\infty} \frac{(-1)^{r-1}re^{-2\pi t r}g(a)}{4t^2r^2+(N+da^2)^2} - \frac{\coth(\pi t)}{\pi}\sum\limits_{r=1}^{\infty} \sum\limits_{a=1}^{\infty} \frac{(-1)^{r-1}e^{-2\pi t r}(4t^2r^2-(N+da^2)^2)g(a)}{(4t^2r^2+(N+da^2)^2)^2} &\\& + \frac{\sinh(\pi t) G_{-N,t,k}}{8\sqrt{k} t } \sum\limits_{a=1}^{\infty} \frac{g(a)}{\cosh(\frac{\pi da^2}{2t})} + \frac{t\sinh(\pi t)G_{-N,t,k}}{2\sqrt{k}\pi}\sum\limits_{r=0}^{\infty} \sum\limits_{a=1}^{\infty}\frac{(-1)^{r+da-1}(2r+1)e^{-\pi t(2r+1)}g(a)}{t^2(2r+1)^2+d^2a^4} &\\&+\frac{\sinh(\pi t)}{8\sqrt{k}t}\sum\limits_{r=1}^{\infty} \sum\limits_{a=1}^{\infty}\frac{(-1)^rG_{r-N,t,k}g(a)}{ \cosh(\frac{\pi (r+da^2)}{2t})} +\frac{\sinh(\pi t)}{8\sqrt{k}t}\sum\limits_{r=1}^{\infty} \sum\limits_{a=1}^{\infty}(-1)^r\frac{G_{-r-N,t,k}g(a)}{ \cosh(\frac{\pi (r-da^2)}{2t})} &\end{flalign*} \begin{flalign*}&+\frac{t\sinh(\pi t)}{2\sqrt{k} \pi}\sum\limits_{r=1}^{\infty} \sum\limits_{m=0}^{\infty}\sum\limits_{a=1}^{\infty}\frac{(-1)^{m+ad-1}(2m+1)e^{-\pi t(2m+1)}G_{r-N,t,k}g(a)}{t^2(2m+1)^2+(r+da^2)^2} &\\&+\frac{t\sinh(\pi t)}{2\sqrt{k}\pi}\sum\limits_{r=1}^{\infty} \sum\limits_{m=0}^{\infty}\sum\limits_{a=1}^{\infty}\frac{(-1)^{m+ad-1}(2m+1)e^{-\pi t(2m+1)}G_{-r-N,t,k}g(a)}{t^2(2m+1)^2+(r-da^2)^2} , \label{analyticStatementForDifferenceofSquares} \tag{5.14}
\end{flalign*}
where $d,k \in \mathbb{N}$. 
\end{proposition}

Here we shall make some remarks. The evaluations through \eqref{analyticStatement1}-\eqref{analyticStatementForDifferenceofSquares}, and the analytic expression for the function $q_{k,s}$ in \eqref{analyticFormulaForqks}, implies that, with a suitable restriction on $g$, one can always arrive at the analytic expressions of the kinds described in \eqref{analyticStatementForSumofSquares} and \eqref{analyticStatementForDifferenceofSquares} for sums of the form $\sum \frac{g(a)}{b^{4s}} $, where the terms are summed over either the positive integer solutions $(a,b)$ of the Diophantine equation $kb^{2s}+y(a)=N$ or the positive integer solutions $(a,b)$ of the Diophantine equation $kb^{2s}-y(a)=N$, where $y(a)$ is an expression depending on $a$ such that the value of $y(a)$ being integral for all natural numbers $a$.

Next, we point out a particular case of \eqref{analyticStatementForDifferenceofSquares}. We let $d=k=1$ in \eqref{analyticStatementForDifferenceofSquares}, and put $N \to 4N$. Since we have $\sum\limits_{\substack{b^2-a^2=4N \\ (a,b) \in \mathbb{N}^2}} \frac{g(a)}{b^4} = \sum\limits_{\substack{(b,a)=(\frac{N}{d}+d,\frac{N}{d}-d)\\ d|N, \hspace{0.05cm} \frac{N}{d}>d}} \frac{g(a)}{b^4}$, we state 
\vspace{0.7cm}
\begin{proposition}Let $g$ be a function such that for all natural numbers $m$, we have $|g(m)| \leq M$ for some real value $M$. And in addition the series ${\sum\limits_{r=1}^{\infty}}^* \frac{g(r)}{z+r^2}$ converges absolutely for all $z$. Then for all natural numbers $N$ and real $t>0$, we have

\begin{flalign*}
    &\sum\limits_{\substack{d|N \\ \frac{N}{d}>d}} \frac{g(\frac{N}{d}-d)}{(\frac{N}{d}+d)^4}  &\\& =  \frac{\pi^2}{3(e^{2\pi t}-1)} \sum\limits_{a=1}^{\infty} \frac{g(a)}{4N+a^2}+ \frac{\pi - \coth(\pi t)}{2\pi}\sum\limits_{a=1}^{\infty} \frac{g(a)}{(4N+a^2)^2}+ \frac{(1-t^2)\coth(\pi t)}{\pi }\sum\limits_{a=1}^{\infty} \frac{(-1)^{a}g(a)}{(4N+a^2)^2} &\\& + \frac{ \pi^3 \coth(\pi t)}{12t }\sum\limits_{a=1}^{\infty} \frac{(-1)^{a-1}g(a)}{\sinh(\frac{\pi (4N+a^2)}{2t})}  + \frac{\coth(\pi t)(t-\frac{1}{t})}{4}\sum\limits_{a=1}^{\infty} \frac{(-1)^{a}g(a)}{(4N+a^2)\sinh(\frac{\pi (4N+a^2)}{2t})}&\\&+ \frac{\pi \coth(\pi t)}{8}\sum\limits_{a=1}^{\infty} \frac{(-1)^{a}g(a)\cosh(\frac{\pi(4N+a^2)}{2t})}{\sinh^2(\frac{\pi (4N+a^2)}{2t})} - \frac{\pi^2\coth(\pi t)}{3} \sum\limits_{r=1}^{\infty}\sum\limits_{a=1}^{\infty} \frac{(-1)^{r-1}e^{-2\pi t r}(4N+a^2)g(a)}{4t^2r^2+(4N+a^2)^2} &\\& - t\coth(\pi t) \sum\limits_{r=1}^{\infty} \sum\limits_{a=1}^{\infty} \frac{(-1)^{r-1}re^{-2\pi t r}g(a)}{4t^2r^2+(4N+a^2)^2} - \frac{\coth(\pi t)}{\pi}\sum\limits_{r=1}^{\infty} \sum\limits_{a=1}^{\infty} \frac{(-1)^{r-1}e^{-2\pi t r}(4t^2r^2-(4N+a^2)^2)g(a)}{(4t^2r^2+(4N+a^2)^2)^2} &\\& + \frac{\sinh(\pi t) G_{-4N,t}}{8 t } \sum\limits_{a=1}^{\infty} \frac{g(a)}{\cosh(\frac{\pi a^2}{2t})} + \frac{t\sinh(\pi t)G_{-4N,t}}{2\pi}\sum\limits_{r=0}^{\infty} \sum\limits_{a=1}^{\infty}\frac{(-1)^{r+a-1}(2r+1)e^{-\pi t(2r+1)}g(a)}{t^2(2r+1)^2+a^4} &\\&+\frac{\sinh(\pi t)}{8t}\sum\limits_{r=1}^{\infty} \sum\limits_{a=1}^{\infty}\frac{(-1)^rG_{r-4N,t}g(a)}{ \cosh(\frac{\pi (r+a^2)}{2t})} +\frac{\sinh(\pi t)}{8t}\sum\limits_{r=1}^{\infty} \sum\limits_{a=1}^{\infty}(-1)^r\frac{G_{-r-4N,t}g(a)}{ \cosh(\frac{\pi (r-a^2)}{2t})} &\\&+\frac{t\sinh(\pi t)}{2 \pi}\sum\limits_{r=1}^{\infty} \sum\limits_{m=0}^{\infty}\sum\limits_{a=1}^{\infty}\frac{(-1)^{m+a-1}(2m+1)e^{-\pi t(2m+1)}G_{r-4N,t}g(a)}{t^2(2m+1)^2+(r+a^2)^2} &\\&+\frac{t\sinh(\pi t)}{2\pi}\sum\limits_{r=1}^{\infty} \sum\limits_{m=0}^{\infty}\sum\limits_{a=1}^{\infty}\frac{(-1)^{m+a-1}(2m+1)e^{-\pi t(2m+1)}G_{-r-4N,t}g(a)}{t^2(2m+1)^2+(r-a^2)^2} . \label{analyticStatementForASumOfDivisors} \tag{5.15}
\end{flalign*} 
\end{proposition}

Finally in the next section, we shall proceed to provide analytic expressions for some particular arithmetic functions.
\vspace{0.7cm}
\section{Some particular arithmetic functions }
The triple sums appearing in \eqref{analyticStatementForSumofSquares} and \eqref{analyticStatementForDifferenceofSquares} can be further simplified for some specific functions $g$, whenever the series $\sum\limits_{a=1}^{\infty} \frac{g(a)}{z^2+(a^2+x)^2}$, $\sum\limits_{a=1}^{\infty} \frac{(-1)^a g(a)}{z^2+(a^2+x)^2}$ can be given in a closed form. We shall now derive \textbf{Proposition 1}.

\begin{proof}[Derivation of Proposition 1]

The desired equality is arrived by means of simple substitution. We put $g(n)=1$ for all natural numbers $n$ in \textbf{Proposition 15} 
 (Notice that $g$ satisfy the conditions described). We separately consider two cases, one for even $d$ and one for odd $d$. Observe that in both cases, the sums  

 \begin{flalign*}
     \sum\limits_{a=1}^{\infty}\frac{1}{t^2(2m+1)^2+(r \pm da^2)^2},\\
     \sum\limits_{a=1}^{\infty}\frac{(-1)^a}{t^2(2m+1)^2+(r \pm da^2)^2},
 \end{flalign*}

 can be given in closed forms by \eqref{Tformula} and \eqref{Vformula}. Substituting these closed forms, we obtain the desired equality.
\end{proof}

 \textbf{Proposition 2} is also arrived in a similar manner.
\begin{proof}[Derivation of Proposition 2]

In the same manner as in the derivation of \textbf{Proposition 1}, the desired equality is arrived by means of simple substitution. We put $g(n)=1$ for all natural numbers $n$ in \textbf{Proposition 16} 
 (Notice that $g$ satisfy the conditions described). We separately consider two cases, one for even $d$ and one for odd $d$. Observe that in both cases, the sums  

 \begin{flalign*}
     \sum\limits_{a=1}^{\infty}\frac{1}{t^2(2m+1)^2+(r \pm da^2)^2},\\
     \sum\limits_{a=1}^{\infty}\frac{(-1)^a}{t^2(2m+1)^2+(r \pm da^2)^2},
 \end{flalign*}

 can be given in closed forms by \eqref{Tformula} and \eqref{Vformula}. Substituting these closed forms, we obtain the desired equality.
\end{proof}

\section*{Declarations}

The authors have no relevant financial or non-financial interests to disclose.

\end{document}